\documentclass[14pt]{amsart}

\usepackage{amssymb}
\usepackage{amsmath, amscd}
\usepackage{graphicx}

\usepackage[colorlinks]{hyperref}

\newtheorem{thm}{Theorem}[section]
\newtheorem{cor}[thm]{Corollary}
\newtheorem{thm-def}[thm]{Theorem-Definition}
\newtheorem{lem}[thm]{Lemma}
\newtheorem{prop}[thm]{Proposition}
\theoremstyle{definition}
\newtheorem{defn}[thm]{Definition}

\theoremstyle{remark}
\newtheorem{rem}[thm]{Remark}
\numberwithin{equation}{section}

\newcommand{\Spec}{{\rm Spec}}

\DeclareMathOperator{\tensor}{\otimes}

\newcommand{\vacB}{|0\rangle}
\newcommand{\vac}{\vacB}

\newcommand{\nc}{\newcommand}

\nc{\Real}{\mathbb R}
\nc{\Cplx}{\mathbb C}
\nc{\ZZ}{\mathbb Z}
\nc{\LL}{\mathbb L}
\nc{\VV}{\mathbb V}
\nc{\MM}{\mathbb M}

\nc{\Field}{\mathbb F}
\nc{\Rat}{\mathbb Q}

\nc{\gr}{{\rm gr\ }}
\nc{\grG}{{\rm gr^{\mathcal{G}}\ }}
\nc{\grF}{{\rm gr^{\mathcal{F}}\ }}
\nc{\grH}{{\rm gr^{\mathcal{H}}\ }}

\nc{\pole}{\Cplx}

\nc{\de}{\partial}
\nc{\oX}{{\mathcal{O}}_{X}}

\nc{\calO}{{\mathcal{O}}}

\nc{\calA}{{\mathcal{A}}}
\nc{\calB}{{\mathcal{B}}}
\nc{\calC}{{\mathcal{C}}}
\nc{\calD}{{\mathcal{D}}}
\nc{\calE}{{\mathcal{E}}}

\nc{\calF}{{\mathcal{F}}}
\nc{\calG}{{\mathcal{G}}}
\nc{\calH}{{\mathcal{H}}}

\nc{\calK}{{\mathcal{K}}}

\nc{\calT}{{\mathcal{T}}}

\nc{\cA}{\mathcal{A}}
\nc{\cB}{\mathcal{B}}
\nc{\cC}{\mathcal{C}}
\nc{\cD}{{\mathcal{D}}}
\nc{\cE}{{\mathcal{E}}}
\nc{\cF}{{\mathcal{F}}}
\nc{\cG}{{\mathcal{G}}}
\nc{\cH}{{\mathcal{H}}}

\nc{\cJ}{{\mathcal{J}}}

\nc{\cL}{{\mathcal{L}}}
\nc{\cM}{{\mathcal{M}}}
\nc{\cN}{\mathcal{N}}
\nc{\cO}{{\mathcal{O}}}
\nc{\cP}{{\mathcal{P}}}
\nc{\cQ}{{\mathcal{Q}}}

\nc{\cS}{{\mathcal{S}}}
\nc{\cT}{{\mathcal{T}}}

\nc{\cV}{{\mathcal{V}}}
\nc{\cZ}{{\mathcal{Z}}}

\nc{\wt}[1]{\Delta_{#1}}

\nc{\op}[1]{{#1}_{(-1)} }

\nc{\opp}[2]{ {#1}_{(  #2  )} }

\nc{\zero}{{}_{(0)} }
\nc{\one}{{}_{(1)} }
\nc{\opm}{{}_{(-1)} }
\nc{\ops}[1]{{}_{({#1})} }

\nc{\Hom}{{\rm Hom}}
\nc{\Aut}{{\rm Aut}}
\nc{\Endom}{{\rm End\,}}
\nc{\End}{{\mathrm End\,}}

\nc{\Iso}{{\rm Iso}}

\DeclareMathOperator{\id}{id}

\nc{\Ind}{{\rm Ind\ }}
\nc{\im}{{\rm Im\ }}

\nc{\spn}{{\rm span}}
\nc{\Der}{{\rm Der\ }}

\nc{\iso}{\simeq}

\nc{\res}{{\rm res}}
\nc{\Resz}{{\rm Res}_z}
\nc{\ta}{{\tilde{a}}}
\nc{\tb}{{\tilde{b}}}
\nc{\ttt}{{\tilde{t}}}
\nc{\ts}{{\tilde{s}}}
\nc{\tg}{{\tilde{g}}}
\nc{\tf}{{\tilde{f}}}

\nc{\txi}{{\tilde{\xi}}}
 \nc{\teta}{{\tilde{\eta}}}

\nc{\tlam}{{\tilde{\lambda}}}
	
\nc{\la}{\lambda}
\nc{\La}{\Lambda}

\nc{\fg}{\frak{g}}
\nc{\fb}{\frak{b}}
\nc{\fn}{\frak{n}}
\nc{\fh}{\frak{h}}
\nc{\fz}{\frak{z}}
\nc{\fzl}{\fz^{\text{reg},\lambda}}

\nc{\fU}{{\frak{U}}}
\nc{\fM}{{\frak{M}}}

\nc{\catleq}{{\cV}ert^{\leq 1}}
\nc{\slt}{\frak{s}\frak{l}_2}
\nc{\affsl}{\widehat{\frak{s}\frak{l}}_2}

\nc{\prline}{\mathbb{P}^1}

\nc{\Vcrit}{V_{\kappa_c}(\slt)}
\nc{\Symb}{{\rm Symb\  }}

\nc{\omcl}{\Omega_X^{1,\,cl}}
\nc{\homcl}{H^1 (X, \omcl)}
\nc{\funomega}{Fun(\homcl)}

\nc{\cdo}{{\cD}_X^{ch}}
\nc{\dl}{\cD^\cL}
\nc{\twcdo}{{\cD}_X^{ch, tw}}
\nc{\twcdoloc}{{\stackrel{\circ}{\cD}}_X^{ch, tw}}
\nc{\twcdolocflag}{{\stackrel{\circ}{\cD}}_{G/B}^{ch, tw}}
\nc{\gtdo}{{\cD}_X^{tw}}
\nc{\cdofun}[1]{\cD_{#1}\otimes\funomega}

\nc{\Pic}{ {\rm Pic}  }
\nc{\cch}{\Cplx_{\chi(z)}}
\nc{\cps}{\Cplx_{\Psi}}
\nc{\dchmod}{\cdo{\textrm{-}}Mod }
\nc{\dchtwmod}{\twcdo{\textrm{-}}Mod }
\nc{\dmod}{\cD_X{\textrm{-}}Mod}
\nc{\dlmod}{\cD^{\cL}_X{\textrm{-}}Mod}
\nc{\dlm}{\cD^{\cL}_X{\textrm{-}}Mod}

\nc{\dlamod}{\gtdo{\textrm{-}}Mod}
\nc{\tensL}{\otimes_{\oX}\!\cL}
\nc{\ttens}{\overline{\otimes}}
\nc{\bM}{\bar{M}}
\nc{\bg}{\bar{g}}
\nc{\btau}{\bar{\tau}}
\nc{\bxi}{\bar{\xi}}
\nc{\ba}{\bar{a}}

\nc{\abs}[1]{\left\vert#1\right\vert}
\nc{\set}[1]{\left\{#1\right\}}

\nc{\eps}{\varepsilon}
\nc{\To}{\longrightarrow}

\nc{\CDOP}{\cD^{ch}_{\prline}}
\nc{\TWCDOP}{\cD^{ch,tw}_{\prline}}
\nc{\twcdop}{\cD^{ch,tw}_{\prline}}

\nc{\basis}{ \set{\tau^i} }
\nc{\basisom}{\{ \omega_i\}}

\nc{\Lie}{{\rm Lie\,}}
\nc{\lie}[1]{ {\, \rm Lie}^{}_{\,#1}\,  }

\nc{\cfdeg}{{\rm conf.deg}}

\nc{\partialtau}[1]{D_{\de_{#1}}}

\nc{\partialom}[1]{D_{\omega_{#1}}}

\nc{\Fields}{{Fields}}
\nc{\ghat}{{\hat{\fg}}}
\nc{\qf}[2]{\langle {#1}, {#2}  \rangle}
\nc{\vlk}{v_{\lambda, k}}
\nc{\uH}{\underline{\rm{H}}}

\nc{\dexi}{\frac{\partial }{\de x_i}}
\nc{\dexj}{\frac{\partial }{\de x_j}}
\nc{\dexp}{\frac{\partial }{\de x_p}}

\nc{\pair}[1]{\langle #1 \rangle}

\nc{\Heis}{H_X^{\pair{ , }}}

\nc{\homeg}{ H^1 (X, \Omega_X^{ [1, 2 \rangle} ) }
\nc{\homegs}{ (\homeg)^*}

\nc{\lakl}{\langle \lambda^*_k\, |\, \lambda_l^* \rangle}
\nc{\lars}{\langle \lambda^*_r\, |\, \lambda_s^* \rangle}

\nc{\lala}[2]{\langle \lambda^*_{#1}\, |\, \lambda^*_{#2} \rangle}
\nc{\ttau}{ {\tilde{\tau}} }
\nc{\tla}{{\tilde{\lambda}^*}}

\nc{\tfg}{ {\tilde{\fg}} }

\nc{\itaul}{{i^{}_{\tau_l} }}

\nc{\itaum}{{i^{}_{\tau_m} }}

\nc{\ixi}{{i^{}_\xi}}
\nc{\ieta}{{i^{}_\eta}}
\nc{\half}{\frac{1}{2} }

\nc{\tPhi}{\tilde{\Phi}}
\nc{\gV}{{\fg(\cV) }}

\nc{\Vext}{{\cV {\calE} xt}^{\pair{,}}}

\nc{\babla}{{\overline{\nabla}}}

\nc{\bA}{{\overline{\cA}}}
\nc{\bQ}{{\overline{\cQ}}}

\nc{\pairgf}[1]{\pair{#1}_{s, \nabla} }

\nc{\VExt}{{\cV {\mathcal{E}} xt_{\cL}^{\pair{, } }   }} 

\nc{\CExt}{   {\mathcal{C} \mathcal{E}}    xt_{\cL}^{\pair{, }}  }

\nc{\disp}{\displaystyle}

\nc{\Gr}{{\mathcal{G}r}}

\nc{\CDO}{{\mathcal{C}\mathcal{D} \cO}}

\nc{\Om}[3]{{ \Omega_{#3}^{#1} \to \Omega_{#3}^{#2, cl}}}

\nc{\dplus}{\dotplus}
\nc{\dminus}{\dotdiv}

\nc{\dar}{\downarrow}
\nc{\upar}{\uparrow}

\nc{\One}{{\mathbf{1}}}
\nc{\dotdiv}{{\stackrel{.}{-} } }

\nc{\bS}{  {\bar{S}}  }
\nc{\ucA}{\underline{\cA}}
\nc{\uC}{\underline{\Cplx}}
\nc{\bl}{$\bullet$ \,}

\title{Classification of transitive vertex algebroids}
\author{Dmytro Chebotarov}

\begin{document}

\begin{abstract}
We present a classification of  transitive vertex algebroids on a smooth variety $X$
carried out in the spirit of Bressler's classification of Courant algebroids.
In particular, we compute the class of the stack of transitive vertex algebroids. 
We define  deformations of sheaves of twisted chiral differential operators introduced in \cite{AChM}
and use the classification result to describe and classify such deformations.
As a particular case, we obtain a localization of Wakimoto modules at non-critical level on flag manifolds.
\end{abstract}

\maketitle

\section{Introduction}

A vertex algebroid is the algebraic structure induced on a subspace $V_0 \oplus  V_1$
of a vertex algebra $V$.
 The  study of vertex algebroids started with  \cite{GMS}
 where
 the sheaves of chiral differential operators (CDO)
  were defined as the enveloping algebras of
  {exact} vertex algebroids.

Algebras of chiral differential operators
 are sheaves of vertex algebras on smooth varieties
 resembling the associative
 algebras of differential operators in some respects.
   One striking difference from the classical prototype
    is that for some manifolds $X$
   no CDO exists; or if there is one, there may be more than one isomorphism class of such sheaves.
   Speaking in technical language,
   sheaves of CDO form a stack, whose groupoid of global sections may
   be empty or
   have more than one connected component.

   In \cite{GMS} the classification of chiral differential operators was obtained; in particular,
   it was established that a global sheaf of CDO exists on $X$
   if and only if
     $ch_2 (\Omega^1_X) = 0$
     where   $ch_2 (\Omega^1_X)$ is
   the second graded piece
    of Chern character of $\Omega^1_X$.

   This result was re-established by Bressler \cite{Bre}
   in a rather unexpected fashion.
   He noticed that the notion of a vertex algebroid is related to a well-known notion in differential geometry, a Courant algebroid: the latter is a quasi-classical limit of the former.
   He obtained a classification of Courant algebroids extending a fixed Lie algebroid
   and rediscovered  the aforementioned obstruction by
   connecting the existence of a CDO on $X$
   with the existence of certain Courant extensions
    of the Atiyah algebra of the sheaf of 1-forms.

    In both these classification problems the obstruction to global existence
    is a class in
    $
    H^2 (X, \Omega^{2} \to \Omega^{3,cl}).
    $
    This is due to a rather remarkable property of these algebroids:
    one can "twist" an algebroid $\cA$ on $U\subset X$
    by a closed 3-form $\alpha$.
    To be more precise, let us denote by $\VExt$
    (resp. $\CExt$ )\marginpar{!}
    the stack of vertex (resp. Courant) algebroid extensions of a given Lie algebroid
    $\cL  \stackrel{\pi}{\to} \cT_X$
    with an invariant pairing $\pair{,}$ on $\ker \pi$
    (cf. section \ref{Extension-defin Subsection}).
    Then the twisting by 3-form action on
     $\VExt$ and $\CExt$
     extends to  
    an action of a certain stack associated to the complex $\Omega_X^2 \to \Omega^{3, cl}_X$
    (cf. \cite{D})
    and it makes each of those a torsor over the latter.
    By standard abstract nonsense, to every such stack $\cS$ there
    corresponds  a class $cl(\cS)\in H^2(X, \Omega^2 \to \Omega^{3, cl})$
    which vanishes precisely when $\cS$ has a global object.
    For example, the obstruction $ch_2(\Omega^1_X)$ above is exactly the class of
    the stack of exact vertex algebroids  on $X$.

    In this article we classify transitive vertex algebroids.
      Since  exact vertex algebroids   classified in \cite{GMS}
    are, in fact, a particular kind of   transitive vertex algebroids
     ( those whose associated Lie algebroid is the tangent sheaf),
      our classification generalizes that  of \cite{GMS}.

      In particular, we
    compute the class of the stack
    $\VExt$.
    Bressler computes the corresponding class for Courant extensions of $X$ \cite{Bre}
    and proves that
    $cl (\CExt) = -  \frac{1}{2} p_1(\cL, \pair{,})$
    where $p_1(\cL, \pair{,})$ is the Pontryagin class
    associated with the pair
    $(\cL, \pair{,})$,
    a generalization of the familiar first Pontryagin class of a vector bundle,
    defined in {\it loc.cit.}.

   Our main result is Theorem     \ref{intro-thm1} below.
   To prove it we  take up the techniques of Baer arithmetic
   developed by Bressler for Courant algebroids
   and use the classification  of both CDO and Courant algebroids.

    \begin{thm}
    \label{intro-thm1}
%
       The class of $\VExt$ in $H^2(X, \Omega^2 \to \Omega^{3, cl})$ equals
       $$
         cl(\VExt) =  ch_2 (\Omega^1_X)  - \half p_1 (\cL, \pair{,})
       $$ 
    \end{thm}

    It is worthwhile to note that it is possible for a manifold $X$ to have no global CDO
     and Courant extensions of a given Lie algebroid $\cL$, but still have a vertex extension of $\cL$.

     \smallskip

We use the  classification result above to study certain  deformations of sheaves of
{\em twisted chiral differential operators} (TCDO) defined in \cite{AChM}.
A TCDO is defined through a procedure
 that, starting with a CDO produces a sheaf
 which has features of both the original CDO and
  and the Bernstein-Beilinson algebra of twisted differential operators (\cite{BB1}). 
 These sheaves have proved useful in representation theory of affine Lie algebras at the critical level. In particular, one has a localization procedure
  for certain classes of $\hat{\fg}$-modules. \cite{AChM}.

 More explicitly,   a sheaf of TCDO on $X$  is a sheaf of  vertex algebras that locally looks like
$\cD^{ch} \tensor H_X$
where
$\cD^{ch}$ is a sheaf of CDO on $X$
and
 $H_X$ is the algebra of differential polynomials on the space
$H^1 (X, \Omega^{1,cl})$
classifying the {twisted differential operators} on $X$.

 When $X$ is a flag variety, $X = G/B_{-}$, the algebra $H_{G/B}$ is isomorphic to
 $\Cplx[\fh^*]$ where $\fh$ is the Cartan subalgebra of $\fg = \Lie G$.
 Moreover, there is an embedding of affine vertex algebra
 $$
 V_{-h^\vee} (\fg) \to \Gamma(G/B, \cD^{ch, tw}_{G/B})
 $$
 which makes
  the space of sections of the TCDO  over big cell
  $
  \Gamma(U_e, \cD^{ch, tw}_X) \iso  \cD^{ch}(U_e) \tensor H_{G/B}
  $
  a $\fg$-module of the critical level, called the {\em Wakimoto module $W_{0, -h^{\vee}}$}.
  \cite{FF1}
  

 The Wakimoto module $W_{0, -h^\vee}$ is a member of the family 
 $$
 W_{0, k} = \cD^{ch}(U_e) \tensor H_{X, k + h^\vee}
 $$
 where $H_{X, \kappa}$ is the Heisenberg  vertex algebra associated with the space
 $\fh$ with a bilinear form equal to $\kappa$ times the normalized Killing form.

 One might ask whether $W_{0,k}$ with non-critical $k$
 admits a localization similar to that of $W_{0, -h^\vee}$.
   We show that such a sheaf indeed exists on any flag manifold
   and is, in fact, a deformation of the TCDO mentioned above:
   there is one such sheaf for each choice of an invariant inner product $\pair{\cdot, \cdot }$ on $\fg$.
   If $X = \prline$, then we prove that this sheaf is a sheaf of $\widehat{\frak{sl}}_2$-modules
   of level $\pair{\cdot, \cdot } + \pair{\cdot, \cdot }_{crit}$ 
   (cf. Corollary \ref{global_sections_of_algebroid-Cor}).

 More generally, we define a {\em deformation of TCDO} on an arbitrary manifold $X$
 to be the vertex enveloping algebra  of certain transitive vertex  algebroid on $X$.
We apply our main classification  result
(cf. Theorem \ref{intro-thm1} above)
to  classify the deformations.

\smallskip

{\em Acknowledgement.}
 The author would like to thank Fyodor Malikov for immense help and guidance throughout the work.

\bigskip

\section{Preliminaries}

We will recall the basic  notions of vertex algebra 
following the exposition of \cite{AChM}.

All vector spaces will be over  $\pole$. 

\subsection{Definitions and examples}
\label{Definitions}
 Let $V$ be a vector space.

A {\em field} on $V$ is a formal series $$a(z) = \sum_{n\in \ZZ}
a_{(n)} z^{-n-1} \in ({\rm End} V)[[z, z^{-1}]]$$ such that for any
$v\in V$ one has $a_{(n)}v = 0$ for sufficiently large $n$.

Let $\Fields (V)$ denote the space of all fields on $V$.

A {\em vertex algebra}  is a vector space $V$ with the following
data:
\begin{itemize}
  \item a linear map $Y: V \to \Fields(V)$,  
    $V\ni a \mapsto a(z) = \sum_{n\in \ZZ} a_{(n)} z^{-n-1}$
  \item a vector $\vac\in V$, called {\em vacuum vector}
  \item a linear operator $\de: V \to V$, called {\em translation operator}
\end{itemize}
that satisfy the following axioms:
\begin{enumerate}
  \item (Translation Covariance)

 $ (\de a)(z) = \de_z a(z)$

  \item (Vacuum)

  $\vac(z) = \id$;

  $a(z)\vac \in V[z]$ and $a_{(-1)}\vac = a$

  \item (Borcherds identity)
   \begin{align}
   \label{Borcherds-identity}
   &\sum\limits_{j\geq 0} {m \choose j} (a\ops{n+j} b )\ops{m+k-j}\\
   =& \sum\limits_{j\geq 0} (-1)^{j} {n \choose j}\{ a\ops{m+n-j} b\ops{k+j} - (-1)^{n}b\ops{n+k-j} a\ops{m+j}
   \}\nonumber
   \end{align}
\end{enumerate}


A vertex algebra $V$ is {\em graded} if  $V = \oplus_{n\geq 0}V_n$
and for $a\in V_i$, $b\in V_j$ we have
$$a_{(k)}b \in V_{i+j - k -1}$$ for all $k\in \ZZ$. (We put $V_i  = 0$ for $i<0$.)

  All vertex algebras in this article will be graded. 

We say that a vector $v\in V_m$ has {\em conformal weight} $m$ and
write $\Delta_v = m$.

If $v\in V_m $ we denote $v_k  = v_{(k - m +1)}$, this is the
so-called conformal weight notation for operators. One has 
$$v_k V_m
\subset V_{m -k}.$$

 A {\em morphism} of vertex algebras is a map $f: V \to W$ that preserves vacuum and satisfies $f(v_{(n)}v') = f(v)_{(n)}f(v')$.


A {\em module} over a vertex algebra $V$ is a vector space $M$
together with a map
\begin{equation}
\label{def-vert-mod-1} Y^M: V \to \Fields(M),\;  a \to Y^M(a,z) =
\sum_{n\in \ZZ} a^M_{(n)}z^{-n-1},
\end{equation}
 that satisfy the following axioms:
\begin{enumerate}
  \item $\vac^M (z)   = \id_M  $
  \item (Borcherds identity)
  \begin{align}\label{def-vert-mod-2}
     &\sum\limits_{j\geq 0} {m \choose j} (a^{}_{(n+j)} b
     )^M_{(m+k-j)}\\
  = &\sum\limits_{j\geq 0} (-1)^{j} {n \choose j}\{ a^M_{(m+n-j)} b^M_{(k+j)} - (-1)^{n}b^M_{(n+k-j)} a^M_{(m+j)} \}\nonumber
  \end{align}
\end{enumerate}

\smallskip

A module $M$ over a graded vertex algebra $V$ is called {\em graded}
if $M = \oplus_{n\geq 0} M_n$ with
 $v_{k}M_l  \subset M_{l-k}$  (assuming $M_{n} = 0$ for negative $n$).

 A {\em morphism of modules} over a vertex algebra $V$ is a map $f: M \to N$
 that satisfies $f(v^M_{(n)}m) = v^N_{(n)}f(m)$ for $v\in V$, $m\in M$.
$f$ is {\em homogeneous} if $f(M_k)\subset N_k$ for all $k$.

\smallskip


\subsubsection{Commutative vertex algebras.}
\label{Commutative_vertex_algebras} 
A vertex algebra is said to be
  {\em commutative} if $a_{(n)} b =0$ for $a$, $b$ in $V$ and $n\geq 0$.
   The structure of a commutative vertex algebras is equivalent to one 
 of commutative associative algebra with a derivation.

If $W$ is a vector space we denote by $H_W$ the algebra of
differential polynomials on $W$. As an associative algebra it is a
polynomial algebra in variables $x_i$, $\de x_i$, $\de^{(2)}x_i$,
$\dots$ where $\set{x_i}$ is a basis of $W^*$. A commutative vertex
algebra structure on $H_W$ is uniquely determined by attaching the
field
 $x(z) = e^{z\de}x_i$  to $x\in W^*$.

 $H_W$ is equipped with grading such that
 \begin{equation}
 \label{weights-0-1-diff-poly}
 (H_W)_0=\pole,\; (H_W)_1=W^*.
 \end{equation}

\subsubsection{Beta--gamma system.}
\label{beta-gamma-system} 
Define the Heisenberg Lie algebra  to be
the algebra with generators $a^i_n$, $b^i_n$, $1\leq i\leq N$ and
$K$ that satisfy $[a^i_m, b^j_n] = \delta_{m, -n} \delta_{i,j} K$,
$[a^i_n, a^j_m] = 0$, $[b^i_n, b^j_m] = 0$.

Its Fock representation $M$ is defined to be the module induced from
the  one-dimensional representation $\Cplx_1$ of its subalgebra
spanned by $a^i_n$, $n\geq 0$, $b^i_m$, $m>0$ and $K$ with $K$
acting as identity and all the other generators acting as zero.

The {\em beta-gamma system} has  $M$ as an underlying vector space,
the
 vertex algebra structure being determined  by
assigning the fields
$$a^i(z) = \sum a^i_n z^{-n-1}, \ \ b^i(z) = \sum b^i_n z^{-n}$$
to $a^i_{-1}1$ and $b^i_01$ resp., where $1\in \Cplx_1$.

This vertex algebra is given a grading so that the degree of
operators $a^i_n$ and $b^i_n$ is $n$. In particular,
\begin{equation}
\label{weight-0-1-b-g} M_0=\pole[b_0^1,...,b_0^N],\;
M_1=\bigoplus_{j=1}^{N}(b^j_{-1}M_0\oplus  a^j_{-1}M_0).
\end{equation}

\subsection{Vertex algebroids}
\label{Vertexalgebroids}

\subsubsection{Definition}
\label{def-of-vert-alg}
	Let  $V$ be a  vertex algebra. 

	Define a 1-truncated vertex algebra to be a sextuple
	 $(V_0\oplus V_1, \vac, \de, \opm, \zero, \one )$ 
	 where the operations $\opm, \zero, \one$ satisfy all the axioms
	   of a vertex algebra that make sense upon restricting to the subspace $V_0 + V_1$.
	    (The precise definition can be found in \cite{GMS}).
	   The category  of 1-truncated vertex algebras will be denoted $\cV ert_{\le 1}$.

The definition of vertex algebroid is a 
reformulation of that of a sheaf of 1-truncated vertex algebras. 

Let $(X, \cO_X)$ be a  space with a sheaf of $\Cplx$-algebras.

A {\em vertex  $\cO_X$-algebroid} is a sheaf $\cA$ of  $\Cplx$-vector spaces equipped with $\Cplx$-linear maps 
$\pi: \cA \to \cT_X$ and $\de: \cO_X \to \cA$ satisfying $\pi \circ \de = 0$
and with    
operations
$
\opm :  \cO_X \times \cA   \To \cA
$,
$
\zero : \,  \cA \times \cA    \To \cA
$,
$
\one :  \,  \cA \times \cA  \To \cO_X
$
satisfying axioms:
\begin{eqnarray}
\label{v-assoc} 
f\opm(g\opm v) - (fg)\opm v & = & \pi(v)(f)\opm \partial(g) +
\pi(v)(g)\opm \partial(f)
\\
\label{zero-minusone}
 x \zero ( f\opm y )  & = & \pi(x)(f)\opm y + f\opm ( x \zero y )
\\
\label{symm-zero}
 x \zero y  +  y \zero x & = & \partial (  x \one y )
\\
\pi(f\opm v) & = & f\pi(v) 
\\
\label{opm-one}
 ( f\opm x) \one y   & = & f ( x\one y)  -
\pi(x)(\pi(y)(f)) 
\\
 \label{zero-one} 
\pi(v)( x \one y ) & = & (  v \zero x) \one y  +
x \one (v \zero y)
 \\
 \label{d-derivation}
\partial(fg) & = & f\opm \partial(g) + g\opm \partial(f) 
 \\
 v \zero \partial(f)  & = & \partial(\pi(v)(f)) 
 \\
 \label{va-1d}
 v \one \partial(f)   & = & \pi(v)(f)
 \end{eqnarray}
for $v,x,y\in\cA$, $f,g\in\cO_X$.
The map $\pi$ is called the {\em anchor} of $\cA$.

If $\cV = \bigoplus_{n \ge 0} \cV_n$ 
is a (graded) sheaf of vertex algebras with $\cV_0 = \cO_X$, 
then  $\cA = \cV_1$ is a  vertex algebroid with $\de$ equal to the translation operator 
and $\pi$ sending $x\in \cV_1$ to the derivation $f \mapsto x \zero f$.

\subsubsection{Associated Lie algebroid} 
Recall that  a {\em Lie algebroid} 
  is a sheaf of $\cO_X$-modules $\cL$ equipped with a Lie algebra bracket $[, ]$
 and a morphism
$
\pi: \cA \to \cT_X
$
of Lie algebra and $\cO_X$-modules
called {\em anchor}
that satisfies
   $
   [x, ay] = a[x,y] + \pi(x) (a) y
   $,
    $x,y \in \cA$, $a\in \cO_X$.

 \smallskip
 
 If $\cA$ is a vertex algebroid, then 
 the operation $\zero$ descends to that on
 $\cL_{\cA}= \cA / \cO_X \opm \de \cO_X$ 
 and makes it into
  a Lie algebroid, with the anchor induced by that of $\cA$.
  $\cL_{\cA}$  is called {\em the associated Lie algebroid of} $\cA$.

\subsubsection{}
  A vertex (resp., Lie) algebroid is {\em transitive},   if its  anchor map $\pi$ is surjective.

  Being a derivation, see (\ref{d-derivation}),
  $\de: \cO_X \to \cA$ lifts to $\Omega^1_X \to \cA$.
 It follows from (\ref{va-1d}) that if $\cA$ is transitive, then
 $\Omega^1_X \iso \cO_X \opm \de \cO_X$ and $\cA$ fits into an exact sequence
 $$
 0 \To \Omega^1_X \To \cA  \To \cL \To 0,
 $$
$\cL = \cL_{\cA}$ being an extension
  $$
  0 \To \fh(\cL)  \To \cL   \To \cT_X  \To 0
  $$
 where $\fh(\cL) := \ker ( \cL \stackrel{\pi}{\To}  \cT_X )$ is an $\cO_X$-Lie algebra. 
 
%
 
\smallskip

 Note that the pairing $\one$ on $\cA$ induces 
 a symmetric   $\cL_{\cA}$-invariant $\cO_X$-bilinear pairing on $\fg(  \cL_{\cA})$
 which will be denoted by $\pair{ , }$.
 
 We regard the pair $(\cL_{\cA}, \pair{, } )$ as "classical data" underlying the vertex algebroid $\cA$.

\subsubsection{Truncation and vertex enveloping algebra functors}

There is an obvious  truncation functor $$t: \cV ert \to \cV
ert_{\leq 1} $$ that assigns to every vertex algebra a 1-truncated
vertex algebra.
  This functor admits a left adjoint \cite{GMS}
  $$ u:  \cV ert_{\leq 1}  \to \cV ert $$
   called a {\em vertex enveloping algebra functor}.


These functors have evident sheaf versions. 
  In particular, one has the functor
\begin{equation}
U: \cV ert\calA lg \To Sh   \cV ert 
\end{equation}
from the category of vertex algebroids to the category of  sheaves of vertex algebras.

\bigskip

\subsection{Courant algebroids}
We give a definition of a Courant algebroid following \cite{Bre}; see also \cite{LWX}.
\label{Cour-alg-section}

A {\em Leibniz algebra} over $k$ is a $k$-vector space $A$ with a bracket $[, ] :A \tensor_k A \to A$ satisfying
$$
[x, [y,z]] = [[x,y],z] + [y, [x,z]].
$$
The bracket is not assumed to be skew-commutative.

\smallskip

A {\em Courant $\cO_X$-algebroid} is an $\cO_X$-module $\cQ$
equipped with
\begin{enumerate}
\item a structure of a Leibniz $\mathbb{C}$-algebra
$
[\ ,\ ] : \cQ\otimes_\mathbb{C}\cQ \to \cQ \ ,
$
\item
an $\cO_X$-linear map of Leibniz algebras (the anchor map)
$
\pi : \cQ \to \cT_X \ ,
$
\item
a symmetric $\cO_X$-bilinear pairing
$
\pair{, } : \cQ\otimes_{\cO_X}\cQ \to \cO_X \ ,
$
\item
a derivation
$
\partial : \cO_X \to \cQ
$
\end{enumerate}
which satisfy
\begin{align}
\pi\circ\partial & =  0
 \label{complex}
 \\
\left[q_1,fq_2\right] & =  f[q_1,q_2] + \pi(q_1)(f)q_2
 \\
\langle [q,q_1],q_2\rangle + \langle q_1,[q,q_2]\rangle & =  \pi(q)(\langle q_1, q_2\rangle)
\\
\left[q,\partial(f)\right] & =  \partial(\pi(q)(f))
\\
\langle q,\partial(f)\rangle & =  \pi(q)(f) \label{ip-o}\\
\left[q_1,q_2\right] + [q_2,q_1] & =  \partial(\langle q_1, q_2\rangle) \label{ip-symm}
\end{align}
for $f\in\cO_X$ and $q,q_1,q_2\in\cQ$.

A morphism of Courant $\cO_X$-algebroids is an $\cO_X$-linear map
of Leibnitz algebras which commutes with the respective anchor
maps and derivations and preserves the respective pairings.

   A {\em connection} on a Courant algebroid $\cQ$ is an $\cO_X$-linear section $\babla$ of the anchor map such that
    $
    \pair{\babla (\xi), \babla (\eta) } = 0
    $.

If $\cQ$ is a Courant algebroid, then $\cL_{\cQ}= \cQ / \cO_X\de \cO_X$ is a Lie algebroid;
 it is called {\em the associated Lie algebroid of} $\cQ$.
 The pairing $\pair{ , }$ on $\cQ$ induces a $\cL_{\cQ}$-invariant pairing on $\fg(\cL_{\cQ})$
 which will be denoted $\pair{ , }$.


%
%
\bigskip

\subsection{The category of vertex extensions}
\label{Extension-defin Subsection}
Let $\cL$ be a transitive Lie algebroid.
A {\em vertex extension} of $\cL$ is a vertex algebroid $\cA$ with an isomorphism of Lie algebroids
$\phi : \, \cL_\cA \to \cL$. 
In what follows we will always identify $\cL_\cA$ and $\cL$ via $\phi$.

A morphism of vertex extensions of $\cL$
 is a morphism of vertex algebroids
 $f: \cA \to \cA'$
 which induces the identity map on $\cL$.
 Thus $f$ fits into a diagram
 $$
 \begin{CD}
 0 @>>>   \Omega^1_X  @>>> \cA   @>>>   \cL   @>>>   0  \\
 @.  @|  @VV f V        @|  @.     \\
 0 @>>>   \Omega^1_X  @>>> \cA'    @>>>   \cL   @>>>   0
 \end{CD}
 $$

Vertex extensions of $\cL$ on $X$ form a category $\cV \cE xt_{\cL} (X)$; clearly, it is a groupoid.

One can consider the category of vertex extensions of $\cL |_U$ on $U$ for any open subset $U \subset X$. 
These categories with the obvious restriction functors
 form a stack on the Zariski topology of $X$, to be denoted $\cV \cE xt_\cL$

\smallskip

Let $\cA$ be a vertex extension of $\cL$.
Denote $\tfg_\cA := \ker (\pi: \cA  \to \cT_X)$;
it is an extension
$$
\begin{CD}
0 @>>> \Omega^1_X @>>> \tfg_\cA @>>> \fg @>>>   0
\end{CD}
$$
 It is easy to see that the operation $\one$ satisfies
 $
 \tfg_\cA  \one \Omega^1 = 0
 $,
 and, therefore, induces a (symmetric,  $\cO_X$-bilinear) pairing
$$
\pair{, }: \, \fg \times \fg \to  \cO_X
$$
If $f: \cA \to \cA'$ is a morphism of extensions, $f$ induces the identity map on $\fg$; it also preserves $\one$. Therefore
$\cA$ and $\cA'$ must have the same pairing $\pair{,}$ on $\fg$. 
It follows that the groupoid
$
\cV \cE xt_{\cL}(X)
$
is a disjoint union   
	$$
	\cV\cE xt_{\cL} (X)  = 
	\coprod\limits_{\pair{,} }
	 \cV \cE xt_{\cL}^{\pair{,}}(X)
	$$ 
where
$
\cV \cE xt_{\cL}^{\pair{,}}(X)
$
is  the full subcategory of vertex extensions of $\cL$ whose induced pairing on $\fg$ is $\pair{,}$.
 Such extensions will be called {\em vertex extensions of $(\cL, \pair{,})$}.





\smallskip

Similarly, we define the notion of a {\em Courant extension} of $\cL$ on $X$ and that of a morphism of Courant extensions,
 the categories $\CExt(U)$ and $\CExt(U)$, $U\subset X$.

\smallskip

\subsection{Chiral differential operators} 
\label{Chiral_differential_operators}
Vertex extensions of $\cT_X$
are called 
{\em exact vertex algebroids}.
Their vertex enveloping algebras, sheaves of {\em chiral differential operators} (CDO). 
 were first introduced in \cite{MSV}
and classified in \cite{GMS}. 
Let us recall the main classification result.

\smallskip


Let us call a smooth affine variety $U= \Spec A $
   {\em suitable for chiralization}
 if $Der(A)$ is a free $A$-module admitting an abelian frame $\{\tau_1,...,\tau_n\}$.
 In this case
  there is a CDO over $U$, which is uniquely determined by the
  condition that
  $(\tau_i)_{(1)}(\tau_j)=(\tau_i)_{(0)}(\tau_j)=0$.
  Denote this CDO by $D^{ch}_{U,\tau}$.

\begin{thm}
\label{class-cdo-local}
 Let $U=\Spec A$ be suitable for
chiralization with a fixed abelian frame $\{\tau_i\}\subset Der A$.

(i) For each closed 3-form $\alpha\in\Omega^{3,cl}_A$ there is a CDO
over $U$ that is uniquely determined by the conditions
\[
(\tau_i)_{(1)}\tau_{j}=0,\;(\tau_i)_{(0)}\tau_{j}=\iota_{\tau_i}\iota_{\tau_j}\alpha.
\]
Denote this CDO by $\cD_{U,\tau}(\alpha)$.

(ii) Each CDO over $U$ is isomorphic to $\cD_{U,\tau}(\alpha)$ for
some $\alpha$.

(iii) $\cD_{U,\tau}(\alpha_1)$ and $\cD_{U,\tau}(\alpha_2)$ are
isomorphic if and only if there is $\beta\in\Omega^{2}_A$ such that
$d\beta=\alpha_1-\alpha_2$. In this case the isomorphism is
determined by the assignment
$\tau_i\mapsto\tau_i+\iota_{\tau_i}\beta$.
\end{thm}

If $A=\pole[x_1,...,x_n]$, one can choose $\partial/\partial x_j$,
$j=1,...,n$, for an abelian frame and check that the beta-gamma
system $M$ of sect. \ref{beta-gamma-system} is a unique up to
isomorphism CDO over $\pole^n$. A passage from $M$ to Theorem
\ref{class-cdo-local} is accomplished by the identifications
$b^j_01=x_j$, $a^j_{-1}1=\partial/\partial x_j$.

\bigskip

\section{Classification of transitive vertex algebroids}
\label{Classif_vertex_algebroids-Section}
In this section we  present a classification of transitive vertex algebroids in the spirit of \cite{Bre}.

First 
we recall the definition of a
 $\Gr(\Omega_X^{[2,3>})$-gerbe given in \cite{GMS};
 one of the results of \cite{Bre}
 is that  $\VExt$ is a  $\Gr(\Omega_X^{[2,3>})$-gerbe. 
 
In section \ref{linear-algebra-section} we describe the core tool of the classification method: 
the ``addition" operation on various algebroids.
It enables us to construct a vertex extension 
starting from a Courant extension and an exact vertex algebroid.
 With this tool in hands we are able to compute
the class of the stack $\VExt$ 
(Theorem \ref{obstrcution-VExt-thm}).


%

\subsection{Gerbes and torsors}
\label{gerbes section}
\subsubsection{Twisting by a 3-form}
\label{twisting3-subsection}
Let $\cA = (\cA, \opm, \zero, \one, \de, \pi)$ 
be a vertex  extension of $\cL$ on $U\subset X$
 and let $\alpha \in  \Omega^{3,cl}(U)$.
Define an operation
$
{}_{(0) + \alpha} : \cA \times \cA  \to \cA
$
by
\begin{equation}
\label{twisting-by-3-form}
x \, {}_{(0) + \alpha} \, y = x \zero y + \iota_{\pi(x)} \iota_{\pi(y)}  \alpha
\end{equation}

\begin{lem}
Let $\alpha \in \Omega^{3, cl}(U)$. Then:

(1)\,
$
\cA \dplus \alpha := 
   (\cA,     \opm,      {}_{(0) + \alpha},  \one, \de,  \pi )
$
is a vertex extension of $\cL$ on $U$.

(2)\,
The assignment
$\cA \mapsto \cA \dplus \alpha$
 can be extended to an auto-equivalence
  \begin{equation}
  \label{functor-dplus-alpha}
  ? \dplus \alpha : \, \VExt \to   \VExt
 \end{equation}
%
\end{lem}

{\em Proof.}
 The proof  of (1) is the same as in the case of cdo (\cite{MSV, GMS}) or Courant algebroids (\cite{Bre}). 
 To see (2), note that every morphism $f: \cA \to \cA'$ is automatically a morphism 
$\cA  \dplus \alpha \to \cA'  \dplus \alpha$; this tautological action on morphisms makes 
 $?\dplus \alpha$ a functor; 
 the composition $(? \dplus (-\alpha)) \circ (? \dplus \alpha)$ is the identity functor of $ \VExt$. $\qed$

\smallskip

It is clear that the functors
$? \dplus \alpha$, $\alpha \in \Omega^{3,cl}(U)$
define an action of the abelian group $\Omega^{3,cl}(U)$ on the category 
$\cV \cE xt_{\cL}(U)$.
Let us show that this action in fact extends to  an action of a category.


\smallskip

 For an open subset $U\subset X$ define a category
 $
 \Gr (\Omega^{[2,3>})(U)
 $
as follows.
The objects of 
  $
 \Gr (\Omega^{[2,3>})(U)
 $
 are elements $a \in \Omega^{3, cl}(U)$; 
 the morphisms
$$
 Hom (\alpha, \alpha') = \set{ \beta \in \Omega^2 (U)  : \  d\beta= \alpha' - \alpha},
 $$
the composition being the addition in $\Omega^{2}(U)$. 

It is clear that
  $
 \Gr (\Omega^{[2,3>})(U)
 $
 is a groupoid. 
The groupoids 
 $
 \Gr (\Omega^{[2,3>})(U)
 $
 form a prestack
  $
 \Gr (\Omega_X^{[2,3>})
 $;
  the addition of 3-forms gives it the structure of a {\em Picard prestack}.
 See \cite{D}, section 1.4 for generalities on Picard stacks.
  

\smallskip

For $f : \cA \to \cA'$ and  $\beta : \alpha \to \alpha'$ define
\begin{equation}
\label{f-dplus-beta-defn}
(f \dplus \beta)(x) = f(x) + \iota_{\pi(x)} \beta
\end{equation}


\begin{prop}

(1)\,
$f \dplus \beta$ is a morphism of vertex extensions
$$
f \dplus \beta : \cA \dplus \alpha  \to \cA' \dplus \alpha'
$$

(2) \, 
The  formulas (\ref{functor-dplus-alpha})
and (\ref{f-dplus-beta-defn})
define a  functor
\begin{equation}
   \nonumber
    \dplus :  \ \cV \cE xt_{\cL}(U) \times  \Gr (\Omega^{[2,3>})(U)  \To \cV \cE xt_{\cL}(U)
\end{equation}
which gives rise to an action of $\Gr (\Omega^{[2,3>})(U) $
on $\cV \cE xt_{\cL}(U)$
\end{prop}

The verification is, again, straightforward and repeats the analogous  discussion in \cite{GMS}.
$\qed$

\bigskip

\subsubsection{$(\Omega^2 \to \Omega^{3,cl})$-gerbes}

\label{gerbe-definition section}

We will say a stack $\mathcal{S}$ over $X$  is a 
{\em 
$\Gr (\Omega_X^{[2,3>})$-gerbe
}
if there is an action  
$
\dplus: \, \mathcal{S} \times  \Gr (\Omega_X^{[2,3>}) \to \mathcal{S}
$
and a cover $\fU = \set{U_i}_{i\in \mathcal{I}}$
such that for any $i\in \mathcal{I}$  and  $x\in \mathcal{S}(U_i)$
 the functor $x \dplus ?: \,  \Gr (\Omega_X^{[2,3>})(U) \to \mathcal{S}(U)$ 
 is an equivalence. 
   (In other words, $\mathcal{S}$ is a torsor over the associated stack).

 \begin{thm} 
 \cite{Bre}
 \label{Thm - VExt and CExt are gerbes}
 The stacks $\VExt$ and $\CExt$, when locally nonempty, are
 $\Gr (\Omega^{[2,3>})$-gerbes.
 \end{thm}


\smallskip

\begin{rem}
The categories 
  $
 \Gr (\Omega^{[2,3>})(U)
 $,
 $U\subset X$
form a  Picard prestack ({\it cf.}\cite{D}, section 1.4.11)
whose associated stack is
the stack of { \em
$(\Omega^2 \to \Omega^{3,cl})$-torsors}.

What Bressler  shows in  \cite{Bre} is
 that this stack is equivalent to the stack 
$\cE \cC \cA_X$
of 
{\em exact Courant algebroids}, and
 that the stacks
 $\VExt$, $\CExt$ 
 are, in fact,
 $\cE\cC\cA_X$-torsors. 
%
\end{rem}


\smallskip

Observe that for $\VExt$ being a $\Gr (\Omega^{[2,3>})$-gerbe
means that for small enough $U\subset X$ and $\cA \in \VExt$
one has an equivalence
$$
\cA \dplus ? : \,   \Gr (\Omega^{[2,3>})(U)  \to  \VExt(U)
$$
In particular, there is an isomorphism
$$
\Hom (\alpha, \alpha') \iso \Hom (\cA + \alpha, \cA + \alpha' )
$$
Under this isomorphism, an element
 $\beta \in \Omega^{2}$ with $d\beta = \alpha' - \alpha$, 
is mapped to the morphism (cf. (\ref{f-dplus-beta-defn}))
\begin{equation}
\label{definition-of-exp}
\exp (\beta):= \id \dplus \beta : \ x \mapsto x + \iota_{\pi(x)} \beta
\end{equation}

The same is true for Courant algebroids and we will use the notation 
$\exp(\beta)$ in both cases.

\bigskip

\subsubsection{The class of a gerbe}

Let $\cS$ be 
a $\Gr (\Omega^{[2,3>})$-gerbe and $\fU$ a cover as in \ref{gerbe-definition section}.
Let us choose an object $x_i \in \cS(U_i)$  for each $i$. 
For each pair $i,j$ we have objects
$x_i |_{U_{ij}}$ and  $x_j |_{U_{ij}}$, and therefore, an isomorphism
\begin{equation}
\eta_{ij} :  x_i |_{U_{ij}} \to x_j |_{U_{ij}} \dplus \alpha_{ij}
\end{equation}
for some $\alpha_{ij} \in \Omega^{3,cl}$. 

The collection  $(x_i, \eta_{ij},  \alpha_{ij})$ is called a {\em trivialization} of $\mathcal{S}$.

We will denote by the same letter $\eta_{ij}$
all of its translates
$$
\eta_{ij} \dplus  \id_\gamma :  x_i |_{U_{ij}} \dplus \gamma
 \to
  x_j |_{U_{ij}} \dplus  (\alpha_{ij} +  \gamma)
$$
for $\gamma \in \Omega^{3,cl}(U_{i})$.

For each triple $i, j, k$ consider the composition (over $U_{ijk} = U_i \cap U_j \cap U_k$) 
$$
\begin{CD}
\eta_{jk} \eta_{ij}  \eta_{ik}^{-1} : \ 
x_k  
@> \eta_{ik}^{-1} >>
    x_i \dplus (-\alpha_{ik}) 
@> \eta_{ij}  >> 
  x_j \dplus (\alpha_{ij} - \alpha_{ik})
 @> \eta_{jk} >>
   x_k \dplus  (\alpha_{ij} + \alpha_{jk} - \alpha_{ik})
\end{CD}
$$
and denote by $\beta_{ijk}$ the element of $\Omega^2(U_{ijk})$ such that
\begin{equation}
\eta_{jk} \eta_{ij}  \eta_{ik}^{-1} = \exp(\beta_{ijk})
\end{equation}
One checks that
\begin{equation}
d_{\check{C}}\beta_{ijk} = 0, \ \ d_{DR} (\beta_{ijk}) = d_{\check{C}} (\alpha_{ij}), \ \ 
d_{DR} (\alpha_{ij}) = 0
\end{equation}
so that the pair 
$
(\alpha_{ij}, \beta_{ijk}) 
$
is an element of $\check{Z}^2 (\fU, \Omega^{2} \to \Omega^{3, cl})$.

By definition,  the {\em class} of  $\mathcal{S}$,
  $cl(\mathcal{S})$,
 is  the class of 
$
(\alpha_{ij}, \beta_{ijk}) 
$
in $H^2(X, \Omega^{2} \to \Omega^{3, cl})$. 
One has the following classical result (cf., e.g., \cite{GMS} for a proof).
\begin{prop}
$\mathcal{S}(X)$ is nonempty if and only if $cl(\mathcal{S}) =0$.$\qed$
\end{prop}

\smallskip

\subsection{The stack $\CExt$}
\label{facts about CExt - section}
As an example, and for future use, 
we recall the construction of a trivialization of the stack $\CExt$ given in \cite{Bre}. 

\smallskip

Let us choose a cover $\fU = \set{U_i}$ such that    $\cT_{U_i}$ 
is free, choose  connections ($\cO_X$-linear sections of the  anchor map)
$$
\nabla_i : \cT_{U_i} \to \cL_{U_i}
$$ 
and identify $\cL_{U_i} \iso \cT_{U_i} \oplus \fg_{U_i}$ via $\nabla_i$.

Define 
$
 c_i = c(\nabla_i) \in \Omega^{2, cl}_{U_i} \tensor_{\cO}  \fg_{U_i}
$
to be the {\em curvature} of the connection $\nabla_i$,     i.e.
$$
c_i (\xi, \eta) = [\nabla_i (\xi), \nabla_i(\eta)] - \nabla_i ([\xi, \eta])
$$

Recall the following 
\begin{thm} \cite{Bre} 
\label{Thm - local existence of CExt}
Let $U \subset X$ and  $\nabla: \cT_U \to \cL_U$
is any connection.

Then
  the category  $\CExt{}(U)$ is nonempty if and only if
  the form
  $  \disp
   \half \pair{c(\nabla) \wedge c(\nabla) }
  $
  is exact.
\end{thm}

\bigskip

Assume that 
$$ 
   \half \pair{c(\nabla_i) \wedge c(\nabla_i) } = dH_i
  $$
for some $H_i \in \Omega^3$.
 Then one can construct a Courant extension $\cQ_{\nabla_i, H_i}$, 
 which is equal to
$ \cL_{U_i} \oplus \Omega^1_{U_i}$ 
as a sheaf of $\cO_U$-modules, and satisfies 
\begin{gather}
\label{Q_NH xi eta}
[\xi, \eta] = [\xi, \eta]_{\cL} + \iota_\xi \iota_\eta H_i, \ \ \ \xi, \eta \in \cT_{U_i}, \\
\label{Q_NH fg pair omega}
<\fg, \Omega^1_U> = <\fg, \, \nabla_i (\cT_U)> =0, \\
\label{Q_NH xi-braket-g}
     [\xi, g]   = [\nabla_i(\xi), g]_{\cL}    -   \pair{\iota_{\xi }c(\nabla_i), g}.
\end{gather}

\smallskip

For each $i, j$ define
$$
A_{ij} = \nabla_i -  \nabla_j \  \in \, \Omega^1_{U_{ij}} \tensor  \fg^{}_{U_{ij}}
$$
\begin{thm} \cite{Bre}
There exists an isomorphism in $\CExt (U_{ij})$
\begin{equation}
\theta_{ij}: \,
\cQ_{\nabla_i, H_i}  
 \stackrel{\sim}{\To} \cQ_{\nabla_j, H_j} \dplus  \alpha_{ij}
\end{equation}
given by
\begin{equation}
\begin{split}
\xi & \mapsto   \xi + A_{ij}(\xi)  - \half \pair{A_{ij}(\xi), A_{ij}  }
\\
g &\mapsto  g - \pair{g, A_{ij} }
\\
\omega&  \mapsto  \omega
\end{split}
\end{equation}
where
\begin{equation}
\alpha_{ij} = \pair{c(\nabla_i) \wedge A_{ij}}  
        -  \half \pair{ [\nabla_i, A_{ij}], A_{ij}} 
             + \frac{1}{6} \pair{ [A_{ij}, A_{ij}], A_{ij}}
             + H_i - H_j
\end{equation}

\end{thm}

\smallskip

The collection $(\cQ_{\nabla_i, H_i},  \theta_{ij}, \alpha_{ij})$ is a trivialization 
of the gerbe $\CExt$.

\bigskip

On triple intersections $U_{ijk}= U_i \cap U_j \cap U_k$ 
the isomorphisms $\theta_{ij}$ satisfy (\cite{Bre})
\begin{equation}
\theta_{jk} \theta_{ij} \theta_{ik}^{-1} = 
  \exp( - <A_{ij} \wedge A_{jk}>) \ \ \ \  \qed
\end{equation} 

\smallskip

Define $\beta_{ijk} = - <A_{ij}  \wedge A_{jk}>$. 

Then
$
 (\alpha_{ij}, \beta_{ijk}) 
$
is a cocycle in 
  $\check{Z}^2(\fU, \Omega_X^2 \to \Omega_X^{3, cl})$.
The corresponding cohomology class
was  
 identified in \cite{Bre} with minus one half of
  {\em the first Pontryagin class $p_1 (\cL, \pair{,})$
   of $(\cL, \pair{,}$
   }.

\begin{thm} 
\cite{Bre}
\label{class-CExt-Theorem}
This class is the class of the stack $\CExt$:
$$
cl(\CExt) = - \half p_1 (\cL, \pair{,}).
$$
\end{thm}    

  
\bigskip

  \subsection{Linear algebra} 
  \label{linear-algebra-section}
   In this section we describe
 the main tool in the proof of the classification result: 
 we define linear algebra-like operations on various algebroids. 
The main technical result to be proved in this section is as follows.

\begin{thm} 
\label{lin-alg-main-thm}
Let $U$ be suitable for chiralization.

Then there exist a functor
$$
\boxplus: \CExt(U) \times \CDO(U)  \To \VExt(U)
$$
 $$
 (\cQ, \cD) \mapsto \cQ \boxplus \cD
 $$
and a functor
$$
\boxminus: \VExt(U) \times \CDO(U)  \To \CExt(U)
$$
 $$
 (\cA, \cD) \mapsto \cA \boxminus \cD
 $$
such that
 for a fixed $\cD \in \CDO(U)$ the functors
	$$
 - \boxminus \cD : \ \VExt(U)  \to \CExt(U) 
     $$
     and
     $$
 - \boxplus \cD : \ \CExt(U)  \to \VExt(U)
        $$
 are  mutually inverse equivalences of $\Gr (\Omega^{[2,3>}(U))$-torsors.
\end{thm}

\smallskip

In fact, the functor $\boxminus$ was defined in \cite{Bre}, together with several versions of $\boxplus$ defined for  various algebroids.  Our $\boxplus$ is just an extension of Bressler's definition.

\bigskip

\subsubsection{Addition}
\label{section-Addition}
Let $\cQ$ be a Courant extension of $\cL$
and $\cD$ a  cdo. 

We describe how to define a vertex extension of $\cL$ 
which can be though of as a ``sum"  of these two structures;
the construction parallels that of the Baer sum of two extensions. 

First, consider the pullback
$
 \bA:= \cQ \times_{\cT} \cD
$
so that a section of $\cA$ is a pair $(q,x)$, $q\in \cQ$, $x \in \cD$ with $\pi(q) = \pi(x)$.

Define  operations
$
\opm: \cO_X \times  \bA \to \bA
$
and
$
\zero, \one : \, \bA \times \bA \to \bA
$
as follows:
\begin{eqnarray}
\label{boxplus-oper-1}
 a\opm (q,x)       & := & (aq, a \opm x)    
 \\
 (q,x) \zero (q', x')   & := & ( \, [ q, q']^{}_\cQ ,  x \zero x'),  
 \\
  (q,x) \one (q', x')   & := & \pair{q, q'}  + x \one x',  
  \\
\pi( (q,x ) )  &:=& \pi(q) = \pi(x), \\
\label{boxplus-oper-5}
\de a   & = & (\de a, 0)
\end{eqnarray}

Note that $\bA$ contains two copies of $\Omega^1$, one from $\cQ$ and the other from $\cD$. 

Let us define $\cQ \boxplus \cD$ to be  the pushout  of $\bA$
with respect to the addition map $+: \Omega^1 \times \Omega^1 \to \Omega^1$
so that one has the following 
$$
\begin{CD}
0 @>>> \Omega^1 \oplus \Omega^1 @>>> \bA @>>> \cL @>>> 0  \\
@.   @VV+V       @VVV       @|         @. \\
0 @>>> \Omega^1 @>>> \cQ \boxplus \cD  @>>> \cL @>>> 0 
\end{CD}
$$

Alternatively, $\cQ \boxplus \cD$ fits into the diagram
$$
\begin{CD}
0 @>>> \tfg \oplus \Omega^1 @>>> \bA @>\pi>> \cT_X @>>> 0  \\
@.   @VV+V       @VVV       @|         @. \\
0 @>>> \tfg @>>> \cQ \boxplus \cD  @>>> \cT_X @>>> 0 
\end{CD}
$$
where the rows are exact and the left square is a push-out square.


\begin{thm}
\label{boxplus-is-valgd-thm}
The operations  (\ref{boxplus-oper-1} - \ref{boxplus-oper-5}) make sense on 
  $\cQ \boxplus \cD$ 
  and give it the structure of  a vertex algebroid
 \end{thm}

{\em Proof.}  
The verification is straightforward. 
As an example, let us show that (\ref{opm-one}) is satisfied.

For $f\in \cO_X$, $q\in \cQ$, $v\in \cD$, one has:
\begin{gather} 
\nonumber
(f \opm (q, v)) \one (q', v')     =   (fq, f\opm v) \one (q', v')  
  = \pair{fq, q'} + (f \opm v) \one v'
  \\
  \nonumber
  = f \pair{q, q'} + f  (v \one v') - \pi(v) \pi(v') (f)
  = f   ((q, v) \one (q', v')) - \pi((q,v)) \pi((q',v')) (f)
\end{gather}
$\qed$


\bigskip

Note that the assignment $(\cQ, \cD) \mapsto \cQ \boxplus \cD$ is naturally a functor
$$
\boxplus: \CExt(U) \times \CDO(U)  \To \VExt(U)
$$
Indeed, let $f \in \Hom_{\calC\calE xt} (\cQ, \cQ')$,
$f \in \Hom_{\CDO} (\cD, \cD')$. In particular, $f$ and $g$ are maps over $\cT$,
so  $(f, g)$ takes $\cQ \times_\cT \cD \subset \cQ \times \cD$ to $\cQ' \times_\cT \cD'$.
Since $f$ and $g$ act as  identity on the subsheaf $\Omega^1$,
$(f,g)$ gives a well-defined map
 between  the pushouts
 $\cQ \boxplus \cD \to \cQ' \boxplus \cD'$
 that will be denoted $ f \boxplus g $.
Finally, 
it remains to note that
the composition is
``coordinate-wise":
 \begin{equation}
 \label{componentwise-composition-for-dplus}
 (f\boxplus g) (f' \boxplus  g') = f f' \boxplus g g'
 \end{equation}
which implies that  $(f,g) \mapsto  f\boxplus g$ is a functor.

\smallskip

Let us note, among the elementary properties of this functor, the following:

\begin{enumerate}
\item
for any $\alpha \in \Omega^{3, cl}_U$,
$\cQ \in \CExt(U)$, $\cD \in \CDO_U$
 one has the equalities
\begin{equation}
\label{boxplus-dplus-compat-eqn}
 (\cQ \boxplus \cD) \dplus \alpha
   \cong
     \cQ \boxplus (\cD \dplus \alpha)
    \cong
    (\cQ \dplus \alpha) \boxplus \cD
\end{equation}
(by definition of $\dplus \alpha$ the three parts of the equation have underlying sheaf $\cQ \boxplus \cD$, one only has to check that the operations coincide).

\smallskip

\item
one has the equality
	\begin{equation}
        \nonumber
	\exp(\beta) \boxplus \id_{\cD} = \exp(\beta ) = \id_{\cQ} \boxplus \exp(\beta)
	\end{equation}
	in
	$\Hom_{\cV \calE xt}( \cQ \boxplus \cD, (\cQ \boxplus \cD) \dplus d\beta)$;
	more generally, 
	\begin{equation}
	    	\label{exp-boxplus-exp}
		\exp(\beta') \boxplus \exp(\beta'') = \exp(\beta' + \beta'')
	\end{equation}
\end{enumerate}

\subsubsection{Subtraction}
\label{section-Subtraction}
 Let  $\cA$ is a vertex extension of $\cL$ and $\cD$ a cdo.
 In  [Bre] it is described how to define a Courant algebroid $\cA  \boxminus  \cD$.
 Let us recall this construction.

Let
$
 \bQ:= \cA \times_{\cT} \cD
$
so that a section of $\bQ$ is a pair $(v,x)$, $v\in \cA$, $x \in \cD$ with $\pi(v) = \pi(x)$.

Define operations
$
\cdot : \cO_X \times  \bQ \to \bQ,
$
$
[,] : \, \bQ \times \bQ \to \bQ,
$
$
\pair{ , }:  \bQ \times \bQ \to \cO_X,
$
$
\pi: \bQ \to \cT,
$
and $\de: \cO_X \to \bQ$
as follows:
\begin{eqnarray}
 a \cdot  (v,x)   & := & (a \opm  v,  \,   a \opm  x )
 \\
 \left[ (v, x) ,  (v', x') \right]  & := & ( v \zero v',  \,  x \zero x')
 \\
  \pair{ (v,x),  (v', x') }   & := & v \one v'   -    x \one x' \\
\pi( (v,x ) )  &:=& \pi(v) = \pi(x) \\
\de a   & = & (\de a, 0)
\end{eqnarray}

Define  $\cA \boxminus \cD$ to be the pushout of $\bQ$
with respect to  the subtraction map $-: \Omega^1 \times \Omega^1 \to \Omega^1$.

One can show that all operations defined above make sense on
 $\cA \boxminus \cD$. One has
\begin{thm} (\cite{Bre}, Lemma 5.6)
The sheaf
$\cA \boxminus \cD$
with the operations defined above
is a Courant algebroid
 \end{thm}



	\smallskip

  \subsubsection{Compatibility of $\boxplus$ and $\boxminus$}	
	
	\begin{thm}
	\label{inverse-equivalences-thm}
	The functors
	$$
 - \boxminus \cD : \ \CExt(U)  \to \VExt(U) 
     $$
     and
     $$
 - \boxplus \cD : \ \CExt(U)  \to \VExt(U)
        $$
        are mutually inverse equivalences of 
        $\Gr (\Omega^{[2,3>}(U))$-torsors.
	\end{thm}

	{\it Proof.} 
	The compatibility of $\boxplus \cD$ and $\boxminus \cD$ with $\Gr (\Omega^{[2,3>}(U))$-action
	follows from properties (\ref{boxplus-dplus-compat-eqn} - \ref{exp-boxplus-exp}) and their obvious analogs for $\boxminus$.
	Let us construct the natural  isomorphisms
	$
	\eta_{\cA}: \, \cA \iso (\cA \boxminus \cD) \boxplus \cD
	$
	where $\cA$ is a vertex extension of $\cL$ and $\cD$ is a cdo.

        Define
	$$
	\eta_{\cA}(v) = ((v, x), x)
	$$
	where $x\in \cD$ is arbitrary.

	To show $\eta_{\cA}$ is well-defined note that for any $x, y \in \cD$ with $\pi(x) = \pi(y) = \pi(v)$ we have $x-y \in \Omega^1$ and
	$$
	((v, x), x) = ((v,  (y-x) + y ),  x)= ((v, y) + (y-x),  x) = ((v, y), x + (y-x)) = ((v, y), y)
	$$
	To verify $\eta_{\cA}$ is a morphism we check
	$$
	a \opm ((v,x),x) = (a(v,x), a \opm x) = ((a\opm v, a\opm x), a\opm x) = \eta_{\cA}(a \opm v)
	$$
	$$
	((v,x),x) \zero ((v',x'), x') = ( [(v,x) , (v', x') ],  x \zero x' ) = ( (v \zero v' , x \zero x' ) ,  x \zero x' ) = \eta_{\cA}(v \zero v')
	$$
	$$
	((v,x),x) \one ((v',x'), x') =\pair{(v,x), (v', x')} + x \one x'  = v \one v'  - x \one x'  + x \one x'  =  v \one v'
	$$
	To check that $\eta_{\cA}$ is an isomorphism, one can check that the map
		$
	\Psi: (\cA\boxminus  \cD)  \boxplus   \cD  \To \cA,
	$
	$
	((v,x),y)  \mapsto v + (y-x).
	$
	is a well-defined inverse to $\Psi$.
	(Note that every section $((v,x),y)$ of $(\cA  \boxminus \cD)  \boxplus  \cD$ can be written as
	$
	((v,x),y) = ((v+(y-x),x+(y-x)),y) = (( v + (y-x), y), y)
	$
	with $v + (y-x)$ independent of the choice of representative $((v,x), y)$).
	
	The construction of the natural isomorphisms
	    $\eta'_\cQ: \, \cQ \to \cQ \boxplus \cD \boxminus \cD$ 
	is analogous and left to the reader. 
	$\qed$

\smallskip

The constructions of sections 
\ref{section-Addition}, \ref{section-Subtraction}
and Theorem  \ref{inverse-equivalences-thm}
 furnish the proof of Theorem \ref{lin-alg-main-thm}.

\bigskip

\subsection{Classification}

 \subsubsection{Local existence}


Let  $U$ be suitable for chiralization and suppose 
$\nabla: \cT_U \to \cL_U$ is a connection. 
\begin{thm}
\label{vert-cour-local-existence}
Then  the following are equivalent:
\begin{enumerate}
  \item The category  $\VExt{} (U)$ is nonempty
  \item  The category  $\CExt{} (U)$ is nonempty
  \item  The  Pontryagin form
  $ \disp
  \half \pair{c(\nabla) \wedge c(\nabla) }
  $
  is exact.
\end{enumerate}
\end{thm}

{\it Proof.} Since $U$ is suitable for chiralization, there exists a CDO $\cD$ on $X$.
Then (1) and (2) are equivalent due to the addition / subtraction operations: 
given a vertex extension $\cA$ there exists a Courant extension $\cQ = \cA \boxminus \cD$  
and vice versa, given $\cQ$ one can produce a vertex extension $\cA = \cQ \boxplus \cD$.
Finally, the equivalence of (2) and (3) is the content of 
Theorem \ref{Thm - local existence of CExt}. 
$\qed$

\bigskip

  \subsubsection{The obstruction}
\label{obstruction-VExt-section}


\begin{thm}
\label{obstrcution-VExt-thm}
Suppose  $\VExt$ is nonempty. Then    its class is equal to
 $$
 cl(\VExt) = 
 -\half p_1(\cL, \pair{,})   +  ch_2 (\Omega^1_X)
 $$
\end{thm}
where 
	$p_1(\cL, \pair{,}) $ is the first Pontryagin class of a Lie algebroid $\cL$ with pairing $\pair{,}$.
	
\bigskip

{\it Proof.} 
What we will be proving is the following:
 $$
 cl(\VExt) = cl(\CExt) + cl( \CDO(X))
 $$
This is indeed sufficient, 
 in view of Theorem \ref{class-CExt-Theorem}
 and the fact that
  $cl( \CDO(X)) = ch_2 (\Omega^1_X)$ \cite{GMS, Bre}

Let $\fU$ be a cover of $X$ by open subsets $U$ suitable for chiralization. 

Since  $\VExt$ is nonempty,  so is $\CExt$.
Suppose we are given a trivialization of the gerbe $\CDO$ and that of $\CExt$.

In other words, we are given a CDO $\cD_i$ and a Courant extension $\cQ_i$ on each $U_i$,
as well as isomorphisms
$$
\eta_{ij} :\  \cD_i |_{U_{ij}}  \To \cD_j |_{U_{ij}} \dplus \alpha^{ch}_{ij}
$$
and
$$
\theta_{ij}: \cQ_i |_{U_{ij}}  \To \cQ_j |_{U_{ij}} \dplus \alpha^{Q}_{ij}
$$
where
$
 \alpha^{Q}_{ij},
\alpha^{ch}_{ij}  \in \Omega^{3, cl}(U_{ij})
$, 
such that
on triple intersections $U_{ijk} = U_i \cap U_j \cap U_k$ one has
$$
\alpha^{ ch}_{ij} + \alpha^{ch }_{jk} = \alpha^{ ch}_{ik}, \ \ \  \ 
\alpha^{\cQ }_{ij} + \alpha^{ \cQ}_{jk} = \alpha^{ \cQ}_{ik}, 
$$
and
\begin{equation}
\label{isos-on-triple}
 \eta_{jk} \eta_{ij}\eta_{ik}^{-1}  = \exp(\beta^{ch}_{ijk}),
\ \ \ \ \
 \theta_{jk} \theta_{ij} \theta_{ik}^{-1}  = \exp(\beta^{Q}_{ijk}),
\end{equation}
for some
$
\beta^{ch}_{ijk} , \beta^{\cQ}_{ijk} \in \Omega^2  (U_{ijk})
$

Then $(\alpha^{ch}_{ij}, \beta^{ch}_{ijk} )$  
and $(\alpha^{\cQ}_{ij}, \beta^{\cQ}_{ijk} )$
are cocycles representing the classes of the gerbes $\CDO_X$ and $\CExt$ respectively.

\smallskip

Now let us construct a trivialization of the gerbe $\VExt$.
Define
$$
\cA_i = \cQ_i \boxplus \cD_i \in \VExt(U_i).
$$
One has the following isomorphisms:
$$
\begin{CD}
\cQ_i  |_{U_{ij}}    \boxplus \cD_i  |_{U_{ij}}
@> \theta_{ij} \boxplus \eta_{ij} >>
 (\cQ_j \dplus  \alpha^{\cQ}_{ij})  |_{U_{ij}}
  \boxplus
   (\cD_j \dplus \alpha^{ch}_{ij})
  | _{U_{ij}}
  @=
  \cQ_j |_{U_{ij}}   \boxplus \cD_j    |_{U_{ij}}   \dplus  ( \alpha^{\cQ}_{ij} + \alpha^{ch}_{ij}))
  \end{CD}
$$
the latter being the identity on the level of vector spaces, by definition of $?\dplus \alpha$
(cf. sect. \ref{twisting3-subsection}).

Thus
$$
\theta_{ij} \boxplus \eta_{ij} :   \cA_i \stackrel{\sim}{\To} \cA_j   \dplus  ( \alpha^{\cQ}_{ij} + \alpha^{ch}_{ij}).
$$

The collection $(\cA_i,  ( \alpha^{\cQ}_{ij} + \alpha^{ch}_{ij}),  \theta_{ij} \boxplus \eta_{ij} )$
 is a trivialization of the gerbe $\VExt$. 
Let  us compute its class.

On triple intersections $U_{ijk} = U_i \cap U_j \cap U_k$ we have
 (cf. (\ref{isos-on-triple}),  (\ref{componentwise-composition-for-dplus}),
  (\ref{exp-boxplus-exp}) )
\begin{equation}
\begin{split}
( \theta_{jk} \boxplus \eta_{jk} )
( \theta_{ij} \boxplus \eta_{ij} )
( \theta_{ik} \boxplus \eta_{ik} )^{-1}
& =
 \theta_{jk}\theta_{ij} \theta_{ik}^{-1}
\boxplus
\eta_{jk}  \eta_{ij}  \eta_{ik}^{-1} \\
& =
 \exp(\beta^{\cQ}_{ijk})
 \boxplus
  \exp(\beta^{ch}_{ijk}) \\
  &=
   \exp(\beta^{\cQ}_{ijk} + \beta^{ch}_{ijk})
\end{split}
\end{equation}
(here, again, we slightly abuse the notation by writing $\theta_{ij}$ for any of its translates under the action of $\Gr(\Omega^{[2,3>})$).

It follows that
$ ( \alpha^{\cQ}_{ij} + \alpha^{ch}_{ij} , \beta^{\cQ}_{ijk} + \beta^{ch}_{ijk} )$
is a cocycle representing the class
of the gerbe $\VExt$.
$\qed$

\bigskip

\section{Deformation of twisted CDO}

\subsection{Twisted chiral differential operators}

In this section we recall the definition of the sheaf $\twcdo$
of twisted chiral differential operators (TCDO) 
corresponding to a given CDO $\cD^{ch}$
on a smooth projective variety $X$.

\subsubsection{The universal Lie algebroid $\cT^{tw}$}
\label{universal-Lie-algd}
The Lie algebroid underlying TCDO is a ``family of all  TDO". 
More precisely, the universal enveloping algebra $\cD_X^{tw}$ of $\cT^{tw}$ possesses the following property:  
for every $\la\in H^1(X,\Om{1}{2}{X})$  there exists an ideal ${\frak{m}}_\la \subset \cD_X^{tw}$
such that the quotient 
$
\cD_X^{tw} / {\frak{m}}_\la
$
is isomorphic to the tdo $\cD_X^{\la}$ corresponding to the class $\la$.

Let us sketch the construction.

\smallskip

Since  $X$ is projective,
$H^1(X,\Om{1}{2}{X})$ is finite-dimensional,
and  there exists an affine cover $\frak{U}$ so that
$\check{H}^1(\frak{U},\Om{1}{2}{X}) =  H^1(X,\Om{1}{2}{X}  )$.

Let $\Lambda =   \check{H}^1(\frak{U},\Om{1}{2}{X})$.
We fix a lifting
$
 \check{H}^1(\frak{U},\Om{1}{2}{X}) \To    \check{Z}^1(\frak{U},\Om{1}{2}{X})
 $
 and identify the former with the subspace of the latter defined by this lifting. 
 Thus, each $\la \in \Lambda$ 
 is a pair of cochains
 $\la = (  (\la^{(1)}_{ij} ),   ( \la^{(2)}_{i}))$ with
 $
  \la^{(1)}_{ij}  \in \Omega^1(U_i \cap U_j),
 $
 $
 \la^{(2)}_{i}  \in \Omega^{2, cl}(U_i),
 $
 satisfying
 $
 d_{DR}  \la^{(1)}_{ij}   =  d_{\check{C}}   \la^{(2)}_{i}
 $
 and
 $
  d_{\check{C}}  \la^1_{ij}   = 0.
 $

For $\la = (\la^{(1)}_{ij},   \la^{(2)}_{i}) \in \Lambda$ denote $\cD^{\la}$
the corresponding sheaf of  twisted differential operators.
One can consider $\cD^{\la}$ as
an enveloping algebra   
of the (Picard) Lie algebroid $\cT^{\la}  = \cD^{\la}_1$ \cite{BB2}.
As an $\cO_X$-module, $\cT^{\la}$ is an extension
$$
0 \To    \cO_X \One     \To    \cT^{\la}     \To \cT_X     \To  0
$$
given by $(\la^{(1)}_{ij})$.
The Lie algebra structure on $\cT^{\la}_{U_i}$ is given by
$
[\xi, \eta]_{\cT^\la}  =  [\xi, \eta]  +  \ixi \ieta \la^2_{i} \One.
$
and $[\One, \cT^{\la}_{U_i}] = 0$.

\smallskip

 Let $\set{  \la^*_i  }$ and $\set{  \la_i  }$ be dual bases 
of $\Lambda^*$ and $\Lambda$ respectively.
Denote by $k$ the dimension of  $\Lambda$.


Define $\cT^{tw}$  to be  an abelian extension
$$
0 \to \cO_X \tensor \Lambda^*     \to   \cT^{tw}_X  \to \cT_X   \to 0
$$
such that $[\Lambda^*, \cT^{tw} ] =0$ and there exist
connections 
$\nabla_i: \cT_{U_i} \to \cT^{tw}_{U_i}$
satisfying
\begin{gather}
\label{univTDO-A}
\nabla_j(\xi) - \nabla_i(\xi) =  \sum_{r} \iota_{\xi} \la^{(1)}_r(U_{ij}) \la^*_r
\\
\label{univTDO-curvature}
[\nabla_i(\xi), \nabla_i(\eta)]- \nabla_i ([\xi,\eta]) =  \sum_r   \iota_{\xi}  \iota_{\eta} \la^{(2)}_r (U_{i})  \la^*_r
\end{gather}

It is clear that the pair 
$(\cT^{tw}, \cO_X \tensor \Lambda^* \hookrightarrow \cT^{tw})$ 
 is independent of the choices made.
 
	We call the universal enveloping algebra 
	$\cD_X^{tw} = U_{\cO_X}(\cT^{tw})$
	\textit{the 	universal sheaf of twisted differential operators}.

\bigskip






\subsubsection{A universal twisted CDO}
\label{A universal twisted CDO}
	Let $ch_2(X)=0$ and fix a CDO $\cD^{ch}_X$. 
	To each such sheaf one
	can attach a {\em universal twisted CDO}, $\cD^{ch, tw}_X$,
	a sheaf of vertex algebras whose "underlying" Lie algebroid is
	$\cT^{tw}_X$.
	 Let us place ourselves in the situation of the previous section, where we had a fixed affine
	cover ${\frak U}=\{U_i\}$ of a projective algebraic manifold $X$,
	 dual bases $\{\lambda_i\}\in
	H^1(X,\Omega^{[1,2>}_X)$, 
	$\{\lambda_i^*\}\in
	H^1(X,\Omega^{[1,2>}_X)^*$, and a lifting
	$H^1(X,\Omega^{[1,2>}_X)\rightarrow Z^1({\frak
	U},\Omega^{[1,2>}_X)$.

        We can assume that $U_i$ are suitable for chiralization. 
         Let us fix, for each $i$, an abelian basis
	$\tau^{(i)}_1,\tau^{(i)}_2,...$ of $\Gamma(U_i,\cT_X)$.
	 Then the CDO $\cD^{ch}$ is given by a
	collection of 3-forms $\alpha^{(i)}\in\Gamma(U_i,\Omega^{3,cl}_X)$
	(cf. sect.\ \ref{Chiral_differential_operators},
	Theorem~\ref{class-cdo-local}) 
	and transition maps
	$
	  g_{ij} : \cD^{ch}_{U_j}|_{U_i \cap U_j } \to  \cD^{ch}_{U_i}|_{U_i \cap U_j }.
	$
	Let us as well fix splittings $\cT_{U_i } \hookrightarrow \cD^{ch}_{U_i}$
	and view $g_{ij}$ 
	as maps 
	$
	 g_{ij}: (\cT_{U_j } \oplus \Omega^1_{U_j} )|_{U_i \cap U_j } \to  (\cT_{U_i } \oplus \Omega^1_{U_i})|_{U_i \cap U_j }
	 $	
	
	The   {\em universal sheaf of twisted chiral differential operators}
	 $\cD^{ch,tw}_X$ corresponding to $\cD^{ch}_X$
	is 
	a vertex envelope of 
	the $\cO_X$-vertex algebroid $\cA^{tw}$ 
	determined by the following:
	\begin{itemize}
	\item
	   $\cA^{tw}$ is a vertex extensions of $(\cT^{tw}_X, 0)$;
	\item
	  there are embeddings  $\cT_{U_i} \hookrightarrow \cA_{U_i}$ such that
	  $$
	   \tau^{(i)}_l \zero \tau^{(i)}_m   
	      =   \iota_{  \tau^{(i)}_l }  \iota_{  \tau^{(i)}_m } \alpha^{ (i) }  
	       + \sum  \iota_{  \tau^{(i)}_l }  \iota_{  \tau^{(i)}_m }   \lambda^{(2)}_k (U_i) \lambda^*_k    
	  $$
	 \item 
	   the transition function  from $U_j$  to $U_i$ is given by
	   \begin{equation}
	     \label{twcdo_trans_fun}
	   g_{ij}^{tw} (\xi)  =  g_{ij}(\xi)  - \sum \iota_{\xi} \lambda_k^{(1)} (U_i \cap U_j) \lambda^*_k 
	   \end{equation}
	\end{itemize}
	See \cite{AChM} for a detailed construction.


%


\subsubsection{Locally trivial twisted CDO}   
\label{Locally trivial twisted CDO}   
   
Observe that there is an embedding
\begin{equation}
\label{emb_coho}
 H^{1}(X,\Omega^{1,cl}_X)\hookrightarrow H^{1}(X,\Omega^{[1,2>}_X)
\end{equation}
The space $H^{1}(X,\Omega^{1,cl}_X)$ classifies {\em locally
trivial} twisted differential operators, those that are locally
isomorphic to $\cD_X$. Thus for each $\lambda\in
H^{1}(X,\Omega^{1,cl}_X)$, there is a unique up to isomorphism TDO
$\stackrel{\circ}{\cD}^{\lambda}_X$ such that for each sufficiently
small open $U\subset X$, $\stackrel{\circ}{\cD}^{\lambda}_X|_U$ is
isomorphic to $\cD_U$. Let us  see what this means at the level of
the universal TDO.

In terms of Cech cocycles the image of embedding (\ref{emb_coho}) is
described by those $(\lambda^{(1)},\lambda^{(2)})$, see section
 \ref{universal-Lie-algd}, where $\lambda^{(2)}=0$, and this
forces $\lambda^{(1)}$ to be closed.  Picking a collection of such
cocycles that represent a basis of $H^{1}(X,\Omega^{1,cl}_X)$ we can
repeat the constructions of sections
 \ref{universal-Lie-algd}
and
\ref{A universal twisted CDO}    
to obtain
   sheaves   
 $\stackrel{\circ}{\cT}^{tw}_X$ and
$\stackrel{\circ}{\cD}^{ch,tw}_X$.
The latter is glued of pieces isomorphic (as vertex algebras) to
$\cD^{ch}_{U_i} \otimes H_X$ with transition functions as in
(\ref{twcdo_trans_fun}); here $H_X$ is the vertex algebra of differential
polynomials on
 $H^{1}(X, \Omega^{1,cl}_X)$. 
 We will call the 
 sheaf
   $\stackrel{\circ}{\cD}^{ch,tw}_X$ the {\em universal locally trivial 
   sheaf of twisted chiral differential
operators.}

  \subsection{TCDO on flag manifolds}

\label{Example:_flag_manifolds.} Let us see what our constructions
give us if $X=\prline$. We have
$\prline=\pole_{0}\cup\pole_{\infty}$, a cover ${\frak
U}=\{\pole_{0},\pole_{\infty}\}$, where $\pole_0$ is $\pole$ with
coordinate $x$, $\pole_{\infty}$ is $\pole$ with coordinate $y$,
with the transition function $x\mapsto 1/y$ over
$\pole^*=\pole_{0}\cap\pole_{\infty}$.

Defined over $\pole_{0}$ and $\pole_{\infty}$ are the standard CDOs,
$\cD^{ch}_{\pole_0}$ and $\cD^{ch}_{\pole_{\infty}}$. The spaces of
global sections of these sheaves are  polynomials in
$\partial^{n}(x)$, $\partial^{n}(\partial_x)$ (or $\partial^{n}(y)$,
$\partial^{n}(\partial_y)$ in the latter case), where $\partial$ is
the translation operator, so that, cf. sect.\
\ref{Chiral_differential_operators},
$$
(\partial_x)_{(0)}x=(\partial_y)_{(0)}y=1.
$$
There is a unique up to isomorphism CDO on $\prline$,
$\cD^{ch}_{\prline}$; it is defined by gluing $\cD^{ch}_{\pole_0}$
and $\cD^{ch}_{\pole_{\infty}}$ over $\pole^*$ as follows
\cite{MSV}:
\begin{equation}         
 \label{CDO_P1_gluing}
x\mapsto 1/y,\; \partial_x\mapsto
(-\partial_{y})_{(-1)}(y^2)-2\partial(x).
\end{equation}
The canonical Lie algebra morphism
\begin{equation}
\label{sl2-p1} sl_2\rightarrow \Gamma(\prline,\cT_{\prline}),
\end{equation}
where
\begin{equation}
 \label{formulas-sl2-p1}
 e\mapsto\partial_x,  \quad   h\mapsto -2x\partial_x,  \quad  f\mapsto -x^2\partial_x,
\end{equation}
$e,h,f$ being the standard generators of $sl_2$, can be lifted to a
vertex algebra morphism
\begin{equation}
\label{verte-sl2-p1} V_{-2}(sl_2)\rightarrow
\Gamma(\prline,\cD^{ch}_{\prline}),
\end{equation}
where
\begin{equation}
\begin{split}
 \label{verte-formulas-sl2-p1}
e_{(-1)}\vac  &\mapsto \,  \partial_x, \\
 h_{(-1)}\vac    & \mapsto  \,  -2(\partial_x)_{(-1)}x, \\
 f_{(-1)}\vac    & \mapsto \,  -(\partial_x)_{(-1)}x^2-2\partial(x).  
\end{split}
\end{equation}
The  twisted version of all of this is as follows (\cite{AChM}).

Since $\dim \prline=1$,
$$
H^1(\prline,\Omega^{1}_{\prline}\rightarrow\Omega^{2,cl}_{\prline})=
H^1(\Omega^{1,cl}_{\prline}),
$$
so all twisted CDO on $\prline$ are locally trivial. Furthermore,
$H^1(\prline,\Omega^{1,cl}_{\prline})=\pole$ and is spanned by the cocycle
$\pole_0\cap\pole_{\infty} \mapsto dx/x$. 
We have
$H_{\prline}=\pole[\lambda^*,\partial(\lambda^*),....]$. Let
$\cD^{ch,tw}_{\pole_0}=\cD^{ch}_{\pole_0}\otimes H_{\prline}$,
$\cD^{ch,tw}_{\pole_{\infty}}=\cD^{ch}_{\pole_{\infty}}\otimes H_{\prline}$ and define $
\cD^{ch,tw}_{\prline}$ by gluing $\cD^{ch,tw}_{\pole_0}$ onto $ \cD^{ch,tw}_{\pole_{\infty}}$
via
\begin{equation}
 \label{gluing_twisted_p1}
\lambda^*\mapsto\lambda^*,\;x\mapsto 1/y,\;\partial_x\mapsto
-(\partial_y)_{(-1)}y^2-2\partial(y) +y_{(-1)}\lambda^*.
\end{equation}
Morphism (\ref{verte-sl2-p1}) ``deforms'' to
\begin{equation}
\label{tw-verte-sl2-p1} V_{-2}(sl_2)\rightarrow \Gamma(\prline, \cD^{ch,tw}_{\prline}),
\end{equation}
\begin{equation}
 \label{tw-verte-formulas-sl2-p1}
e_{(-1)}\vac\mapsto\partial_x, h{(-1)}\vac\mapsto
-2(\partial_x)_{(-1)}x+\lambda^*, f_{(-1)}\vac\mapsto
-(\partial_x)_{(-1)}x^2-2\partial(x)+x_{(-1)}\lambda^*.
\end{equation}
Furthermore, consider
$T=e_{(-1)}f_{(-1)}+f_{(-1)}e_{(-1)}+1/2h_{(-1)}h\in V_{-2}(sl_2)$.
It is known that $T\in\fz(V_{-2}(sl_2))$, the center of
$V_{-2}(sl_2)$, and in fact, the center $\fz(V_{-2}(sl_2))$ equals
the commutative vertex algebra of differential polynomials in $T$.
The formulas above show
\begin{equation}
\label{image-sugawara}
T\mapsto\frac{1}{2}\lambda^*_{(-1)}\lambda^*-\partial(\lambda^*)\in
H_{\prline}.
\end{equation}

All of the above is easily verified by direct computations, cf. \cite{MSV}. The higher rank
analogue is less explicit but valid nevertheless.

\bigskip

Let $G$ be a simple complex Lie group, $B\subset G$ a Borel subgroup, $X=G/B$, the flag
manifold,  $\fg=\text{Lie\;}G$ the corresponding Lie algebra, $\fh$ a Cartan subalgebra.  One
has a sequence of maps
\begin{equation}
 \fh^* \rightarrow H^1(X,\Omega^{1,cl}_X)\rightarrow
H^1(X,\Omega^{1}_X\rightarrow\Omega^{2,cl}_X).
\end{equation}
The leftmost map    attaches to an integral weight $\lambda\in P\subset\fh^*$ the Chern class
of the $G$-equivariant line bundle $\cL_\lambda=G\times_{B}\pole_{\lambda}$, and then extends
thus defined map $P\rightarrow H^1(X,\Omega^{1,cl}_X)$ to $\fh^*$ by linearity. The rightmost
one is engendered by the standard spectral sequence converging to  hypercohomology. It is easy
to verify that both these maps are isomorphisms. Therefore,
\begin{equation}
\label{on-flag-trivial} \fh^* \stackrel{\sim}{\rightarrow}H^1(X,\Omega^{1,cl}_X)
\stackrel{\sim}{\rightarrow}H^1(X,\Omega^{1}_X\rightarrow\Omega^{2,cl}_X),
\end{equation}
and each twisted CDO on $X$ is locally trivial.


Constructed in \cite{MSV} -- or rather in \cite{FF1}, see also \cite{F1} and \cite{GMSII} for
an alternative approach -- is a vertex algebra morphism
\begin{equation}
 \label{verte-g-f}
V_{-h^{\vee}}(\fg)\rightarrow\Gamma(X,\cD^{ch}_X).
\end{equation}
Furthermore, it is an important result of Feigin and Frenkel \cite{FF2}, see also an excellent
presentation in \cite{F1}, that $V_{-h^{\vee}}(\fg)$ possesses a non-trivial center,
$\fz(V_{-h^{\vee}}(\fg))$, which, as a vertex algebra, isomorphic to the algebra of
differential polynomials in $\text{rk}\fg$ variables.
\begin{lem} \cite{AChM}
 \label{tw-verte-g-f}
Morphism (\ref{verte-g-f}) ``deforms'' to
$$
\rho: \;V_{-h\check{ }}(\fg)\rightarrow\Gamma(X,\cD^{ch,tw}_X).
$$
	Moreover, 
	$
	\rho(\fz(V_{-h^{\vee}}(\fg)))\subset H_{X}.
	$
\end{lem}

\smallskip

\subsection{A deformation}

\subsubsection{Motivation: Wakimoto modules}
Let
$X =G/ B_{-}$ be a flag variety
and  $U = NB_{-} \subset X$ 
the {big cell} of $X$.

In virtue of Lemma  \ref{tw-verte-g-f},
the sections $\Gamma(U,   \twcdo)$ become a $V_{-h^\vee}(\fg)$-module, hence a $\fg$-module at the critical level. 
Following \cite{FF1, F3}, we call $\Gamma(U,   \twcdo)$ 
a {\em Wakimoto module of highest weight} $(0, -h^{\vee})$, 
to be denoted $W_{0, -h^\vee}$.

By construction
$
W_{0, -h^\vee} 
= 
\cD^{ch}(U) \tensor H_X
$.
In fact, Feigin and Frenkel proved 
\cite{FF1}
that there exists a whole family of $\fg$-modules
$$
W_{0, k -h^\vee}  = \cD^{ch}(U) \tensor H_k
$$
where $H_k$ is the Heisenberg vertex algebra associated to the space $\fh$
 with bilinear pairing $k\pair{, }_0$, i.e., $k$ times the canonically normalized Killing form. 
The $\fg$-module structure is defined by a vertex algebra morphism
$$
V_{ k -h^\vee}(\fg) \to   \cD^{ch}(U) \tensor H_k
$$
and thus, $W_{0, k -h^\vee}$ is a $\fg$-module of level $k-h^\vee$.

\smallskip

When the level is critical, 
$W_{0, -h^\vee} = \Gamma(U_e, \cD^{ch, tw}_{X})$
One might ask whether  sheaves with an analogous property exist for Wakimoto modules at a non-critical level. 
To be more precise, we are interested in a sheaf $\cV$ of vertex algebras such that:
\begin{itemize}
\item
its sections on the big cell $U$ and its $W$-translates
are isomorphic to the tensor product of vertex algebras $\cD^{ch}({\mathbb{A}}^{\dim \fg /\fb})\tensor H_{k}$, for nonzero $k$;
\item
the associated Lie algebroid of  $\cV$ is the universal tdo $\cT^{tw}$.
\end{itemize}
In other words, $\cV$ is a vertex extension of the pair $(\cT_{G/B}^{tw},  k\pair{,}_{0} )$

We show that such sheaves do indeed exist on $G/B$; 
moreover, the construction is rather general and can be carried out for any variety. 
We call the obtained sheaves the {\em deformations of TCDO} or {\em deformed TCDO};
deformations because they depend on $\pair{,}$ as a parameter, with $\pair{, } = 0$ corresponding to a TCDO.

\subsubsection{Definition}
The discussion above suggests the following definition.

Let $X$ be a smooth projective variety and $\cT^{tw}$ 
the Lie algebroid underlying the universal TDO (cf. section \ref{universal-Lie-algd}).
Recall that $\cT^{tw}_X$ fits into an exact sequence
$$
0 \to \cO_X \tensor \Lambda^* \to \cT^{tw}_X  \to \cT_X  \to 0
$$
where 
$
\Lambda = H^1 (X, \Omega^{1} \to \Omega^{2,cl}).
$

Let us fix a symmetric bilinear pairing 
$\pair{, }: \Lambda^* \times \Lambda^* \to \Cplx$  and extend $\cO_X$-linearly to $\cO_X \Lambda^*$.
\begin{defn} 
We will say that a sheaf $\cV$ is a {\em $\pair{,}$-deformation of TCDO} if $\cV$ is a vertex extension 
of the pair $(\cT^{tw}_X, \pair{, })$.
\end{defn}
Without specifying $\pair{,}$, a {deformation of TCDO} is just a vertex extension 
of the Lie algebroid $\cT^{tw}_X$.

Being vertex extensions, $\pair{,}$-deformations form a stack,  to be denoted 
$$
\cT \CDO_X^{\pair{,}} := \cV\cE xt^{\pair{,}}_{\cT^{tw}} 
$$

\smallskip

  \subsection{Classification of deformations}
 
 We apply the results of sections \ref{obstruction-VExt-section}.
 
 Theorem \ref{obstrcution-VExt-thm} implies that, when 
 $\cT \CDO_X^{\pair{,}}$ is locally nonempty, its class is equal to
 \begin{equation}
 \nonumber
 cl(\cT \CDO_X^{\pair{,}}) = cl(\cC\cE xt_{\cT^{tw}}^{\pair{,}} ) + ch_2 (\Omega^1_X)
 \end{equation}
 
 We are going to use the description of
  $cl(\cC\cE xt^{\pair{,}}_{\cT^{tw}})$
 given in  section \ref{facts about CExt - section}.

  Let us work in the setup of sections \ref{universal-Lie-algd}, \ref{A universal twisted CDO}.
  Thus, we pick a basis 
  $\set{\la_r}$ of $H^1(X, \Omega^{1}_X \to \Omega^{2, cl}_X)$,
   a dual basis 
  $\set{\la^*_r}$ in $H^1(X, \Omega^{1}_X \to \Omega^{2, cl}_X)$,
 and  a lifting
  $H^1(X, \Omega^{1}_X \to \Omega^{2, cl}_X) \to 
      \check{Z}^1(X, \Omega^{1}_X \to \Omega^{2, cl}_X)$,
    so that  each $\la_r$ is a pair of cochains 
    $
   (\la_r^{(1)}, \la^{(2)}_r ) \in \prod \Omega^1(U_{ij}) \times \prod \Omega^{2,cl}({U_{i}})$.
      
      By construction, the Lie algebroid $\cT_X^{tw}$ admits connections
      $\nabla_i : \cT_{U_i} \to \cT^{tw}_{U_i}$ such that
      \begin{equation}
       \label{Aij-for-Ttw}
       A_{ij} := \nabla_i - \nabla_j = - \la^*_k \la^{(1)}_k (U_{ij})
      \end{equation}
      (summation over repeated indices is assumed) and 
      \begin{equation}
      \label{cNabla-for-Ttw}
        c(\nabla_i) =  - \la^*_k \la^{(2)}_k (U_i)
      \end{equation}

\begin{thm}
\label{class-DTCDO-thm}
Let $\pair{} \neq 0$ be a symmetric bilinear form on $\cO_X \tensor \Lambda^*$. Then:

(1) \, 
 $\pair{}$-deformations exist locally on $X$
 if and only if the 4-form
\begin{equation}
\label{def-Pontr-form}
\pair{\la^*_r, \la^*_s} \la^{(2)}_r(U_i) \wedge \la^{(2)}_s (U_i) 
\end{equation}
is exact;

(2)
\, 
Assume (1) and pick, for every $i$, a 3-form $H_i$ such that 
$2dH_i = \pair{\la^*_r, \la^*_s} \la^{(2)}_r(U_i) \wedge \la^{(2)}_s (U_i)$.
Denote
$$
\alpha_{ij}
    =  \half   \pair{\la^*_r, \la^*_s}  \left( 
             \la^{(2)}_r (U_i) 
             +  \la^{(2)}_r  (U_j)
               \right)      
             \wedge \la^{(1)}_s    (U_{ij})         
             + H_i -  H_j
$$
and
$$
  \beta_{ijk}  =    \pair{\la^*_r, \la^*_s}  \la^{(1)}_r(U_{ij})  \wedge  \la^{(1)}_s (U_{jk})
$$
%
  
  Then
   a global $\pair{,}$-deformation exists if and only if the class of the cocycle 
$
(\alpha_{ij},  \beta_{ijk})
$
 in $H^2(X, \Om{2}{3}{X})$
 is equal to $-ch_2(\Omega_X)$
 (minus second graded piece of Chern character of $\Omega^1_X$).
\end{thm}

{\em Proof.}
(1) 
Follows from Theorem \ref{vert-cour-local-existence}, since 
the 4-form (\ref{def-Pontr-form}) is just  the Pontryagin form
$
\half \pair{ c(\nabla_i) \wedge c(\nabla_i)}
$ 
for the Lie algebroid $\cT^{tw}_{X}$.

(2) Using the connections $\nabla_i$  
(and formulas  (\ref{Aij-for-Ttw}), (\ref{cNabla-for-Ttw}))
in the construction of the
section \ref{facts about CExt - section}  one verifies that
the cocycle $(\alpha_{ij},  \beta_{ijk})$ represents the class of $\cC\cE xt^{\pair{,}}_{\cT^{tw}}$.
The statement follows immediately from Theorem    \ref{obstrcution-VExt-thm}
and the fact that 
$
cl(\CDO) = ch_2(\Omega^1_X)
$ \cite{Bre}.
 $\qed$

\bigskip

  \begin{rem}
In the presence of CDO, the classification problem for deformed TCDO becomes one for
Courant extensions of $(\cT^{tw}_X, \pair{,})$, as any CDO $\cD^{ch}$ defines an equivalence
of stacks over $X$
$$
? \boxplus \cD^{ch} : \,  \cC\cE xt^{\pair{,}}_{\cT^{tw}} \to \cT \CDO^{\pair{,}}. 
$$
\end{rem}

\smallskip

\subsection{Deformations of locally trivial TCDO} 
\label{Deformations_of_locally_trivial_TCDO}
Recall from section 
\ref{Locally trivial twisted CDO}
that locally trivial TCDO are constructed in the same way as TCDO by consistently replacing
$H^1(X, \Omega^1 \to \Omega^{2, cl})$
with
$H^1 (X, \Omega^{1, cl})$. 
In particular we construct a Lie algebroid 
$\stackrel{\circ}{\cT}^{tw}$.

We define the corresponding versions of deformations as follows.
A {\em locally trivial deformed TCDO}
is a vertex extension of $\stackrel{\circ}{\cT}^{tw}$.
A {\em locally trivial $\pair{,}$-deformation of TCDO}
is a vertex extension of
	$
	( \stackrel{\circ}{\cT}^{tw}, \pair{,}).
	$

The locally trivial $\pair{,}$-deformations form a stack
$\cT \CDO^{\pair{,}, lt}$.

\smallskip

 Theorem 
\ref{class-DTCDO-thm}
has the following analogue in the locally trivial case:
\begin{thm}      
\label{class-DTCDO-LT-thm}
Let $\pair{} \neq 0$ be a symmetric $\cT$-invariant
 bilinear form on $\cO_X \tensor \Lambda^*$. 
 Then:

(1) \,  $\pair{}$-deformations exist locally on $X$.

(2) \,  
every $\pair{,}$-deformation
 $\cA^{tw, lt}_{\pair{,}}$ 
 is locally isomorphic to
$\cD_U^{ch} \tensor H_{\pair{,}}$ where 
$\cD_U^{ch}$ is a CDO  and 
$H_{\pair{,}}$
is a Heisenberg vertex algebra associated to the space $H^1(X, \Omega^{1,cl})^*$ with 
the bilinear form $\pair{,}$.

(3) \, 
Denote
$$
  \beta_{ijk}  =    \pair{\la^*_r, \la^*_s}  \la^1_r(U_{ij})  \wedge  \la^1_s (U_{jk})
$$
and let $[ (0, (\beta_{ijk})) ]$ stand for the class of $ (0, (\beta_{ijk}))$ 
 in $H^2 (\Omega^2 \to \Omega^{3,cl})$.
 
Then the   class of  
$
\cT \CDO^{\pair{,}, lt}
$
 in $H^2 (\Omega^2 \to \Omega^{3,cl})$
 is given by
$$
cl(\cT \CDO^{\pair{,}, lt})  = ch_2(\Omega^1_X) + [(0, (\beta_{ijk})) ]
$$

\end{thm}

{\em Proof.}
(1) 
By construction, the Lie algebroid $\stackrel{\circ}{\cT}^{tw}$
admits {\em flat} connections $\nabla_i : \cT_{U_i} \to \stackrel{\circ}{\cT}_X^{tw} |_{U_i}$, 
which implies
$
\pair{c(\nabla_i) \wedge c(\nabla_i)} =0.
$
The local existence now follows from Theorem \ref{vert-cour-local-existence}.

(2)
Suppose $\nabla$ is a flat connection on an open set $U \subset X$,
and let $\cQ =\cQ_{\nabla, H}$ be a Courant extension of 
$\stackrel{\circ}{\cT}^{tw}$
 over $U$ 
(cf. \ref{facts about CExt - section}).
 
Then $\cQ \iso \cT_U \oplus  (\cO_U \tensor H^1(X, \Omega^{1,cl})) \oplus \Omega^1_U$
and since $c(\nabla) = 0$ one immediately observes
from (\ref{Q_NH fg pair omega}) 
and (\ref{Q_NH xi-braket-g})
 that the
 constant subsheaf  
    $H^1(X, \Omega^{1,cl})^*$
    ``decouples".  
     It is clear from the construction, that it stays decoupled in $\cQ \boxplus \cD$, for any cdo $\cD$ on $U$. 
     It has a structure of a Courant (equivalently, vertex) algebroid over $\Spec (\Cplx)$ 
 whose vertex envelope is the algebra  $H_{\pair{,}}$.

(3)  
The proof is identical to that of Theorem \ref{class-DTCDO-thm}, Part (2). $\qed$

\bigskip
  
  \subsection{Deformed TCDO on $\prline$}

This is a continuation of Example \ref{Example:_flag_manifolds.}.

Recall that we are  using standard coordinate charts $U_0$ and $U_1$ so that $\prline = U_0 \cup U_1$ with $0\in U_0$, $\infty \in U_1$ and coordinate functions
 $x: U_0 \to \Cplx$ and $y: U_1 \to \Cplx$ 
 with $x = \frac{1}{y}$. 
Denote
$$
\la = 
  \frac{dy}{y} = - \frac{dx}{x}
$$
a cocycle representative of a generator of 1-dimensional  $H^1(\prline, \Omega^{1,cl})$.

By definition,
\begin{equation}
\begin{split}
\cT^{tw}_{U_i} = \cT_{U_i} \oplus \cO_{U_i} \la^* , \ \ i =0,1,  \\
\end{split}
\end{equation}
with Lie bracket defined by
$[\xi, \eta]_\cL = [\xi, \eta]$,  \ 
$[\xi, a  \la^*] = \xi(a) \la^*$.

Let 
$\nabla_i: \cT_{U_i} \to \cT^{tw}_{U_i}$, $i =0,1$ 
be the canonical inclusions. 
The formula (\ref{Aij-for-Ttw}) in this case reads  as
\begin{equation}
\label{nabla_1 - nabla_0}
\nabla_1 - \nabla_0 = \frac{dy}{y} \la^*,
\end{equation}
 which dictates the following gluing map
 $
   g_{01}:  \cT^{tw}_{1}|_{\Cplx^*}   \to \cT^{tw}_{0}|_{\Cplx^*}    
   $
\begin{align}  
     \xi &  \mapsto \xi + \ixi \la \cdot \la^*  \\
     \nonumber
     \la^* &  \mapsto \la^*
\end{align}

In the chosen coordinates, it is
$
\de_y  = -x^2 \de_x + x \la^*.
$

\bigskip

\subsubsection{The deformed TCDO}
We wish to construct a vertex extension of 
$(\cT^{tw}_{\prline}, \pair)$, where $\pair{,}$ is a symmetric 
$\cT^{tw}$-invariant $\cO$-bilinear pairing on 
$\fg(\cT^{tw}) = \cO_X \tensor H^1(X, \Omega^1 \to \Omega^2)^* = \cO \cdot \la^*$.
In this case it is determined by  a number $k\in \Cplx$ assigned to $\pair{\la^* | \la^*}$. 
Let us fix $k$ and  assume $k \neq 0$ ($k=0$ corresponds to the usual TCDO).

Since $\dim \prline =1$, $\Omega^{i}=0$ for $i>1$, in particular 
$H^i (\prline, \Omega^2 \to \Omega^{3,cl}) = 0$ for all $i$.
Therefore there exists a unique vertex extension for  any  pair  $(\cL, \pair{,})$. 
Let $\cA^{tw}_{\pair{,}}$ denote the vertex extension of $(\cT^{tw}_{\prline}, \pair{,})$. 

Denote by $H_{\prline}^{\pair{,}}$ the Heisenberg vertex algebra generated by a filed $\la^*$ satisfying
$
\la^* \one \la^* = \pair{\la^*, \la^*},  \ \ \la^* {}_{(n)} \la^* =0 , n\ne 1.
$
Theorem \ref{class-DTCDO-LT-thm} describes $\cA^{tw}_{\pair{,}}$ locally: 
one has isomorphisms of vertex algebras
$
(\cA^{tw}_{\pair{,}})_{U_i}  \iso \cD^{ch}_{U_i}  \tensor H_{\prline}^{\pair{,}}, \, i=0,1.
$
Some global information is provided by the following
\begin{thm}
(1)\
There are isomorphisms
$\phi_i : \cA^{tw}_{\pair{,}}  |_{U_i} \to  \cD^{ch}_{U_i}  \tensor H_{\prline}^{\pair{,}}$, $i=0,1$,
such that
\begin{gather}
\label{P1-transition-vectorfields-inTheorem}
\phi_0 \phi_1^{-1} (\de_y ) =  - x^2 \de_x  - 2 dx +   x \la^*   +  \frac{1}{2}  \pair{\la^*, \la^*}  dx \\
\label{P1-transition-lambda-inTheorem}
\phi_0 \phi_1^{-1}  (\la^*)  =  \la^* - \pair{\la^*, \la^*} x^{-1} dx
\end{gather}

(2) \, The anchor map of $\cA^{tw}_{\pair{,}}$
induces  a vector space isomorphism
$$
H^0(\prline, \cA^{tw}_{\pair{,}}) \iso H^0 (  \prline, \cT_{\prline}).
$$
\end{thm}

\bigskip

{\em Proof.} (1)
The construction of section \ref{section-Addition}
and the results of section \ref{facts about CExt - section} imply
 that the most general  gluing formula is as follows:
\begin{gather}
 \xi  \mapsto  g_{ij}(\xi ) + A(\xi)  - \frac{1}{2} \pair{ A(\xi), A}  + \iota_{\xi} \beta
 \\
g \mapsto   g - \pair{g, A}
\end{gather}
where 
 $g_{ij}$ is a transition function for a CDO,   
 $\beta \in \Omega^{2}_{U_i \cap U_j}$,
 $A = \nabla_j - \nabla_i$, 
 for some connectoions
   $\nabla_i : \cT_{U_i} \to \cL_{U_i}$.

Applying to our case
and  using  
(\ref{CDO_P1_gluing})  and  (\ref{nabla_1 - nabla_0}), 
 we see that
 \begin{gather}
 \de_y  \mapsto  - x^2 \de_x  - 2 dx +   x \la^*   +  \frac{1}{2}  \pair{\la^*, \la^*}  dx
\end{gather}
and the map
 $
 \tilde{\fg} |_{U_{1}}  \to \tilde{\fg} |_{U_{0}}
 $
 is given by
\begin{equation}
\label{P1-kernel-transition}
\la^* \mapsto   \la^* - \pair{\la^*, \frac{dy}{y}\la^*} = \la^* + k \frac{dx}{x}
\end{equation}

\bigskip

(2) 
The gluing formula 
(\ref{P1-kernel-transition})
 implies that the map $H^0 (\prline, \fg) \to H^1 (\prline, \Omega^1_{\prline})$
 in the long exact sequence associated to
 $0 \To \Omega^1_{\prline} \To \tilde{\fg} \To \fg \To 0$
  is  
 an isomorphism.
 
  Since $H^j(\prline, \Omega^1) =H^k(\prline, \fg)=0$ for $j \ne 1$, $k \ne 0$,
one can conclude that
$
H^i (\prline, \tilde{\fg}) = 0
$
for all $i$.

In turn, the long cohomology sequence associated to the sequence
$$
0 \To  \tilde{\fg}  \To \cA^{tw}_{\pair{,}} \To \cT_{\prline}   \To 0
$$
shows that
$
H^i(\prline, \cA^{tw}_{\pair{,}}) \iso H^i (  \prline, \cT)
$.
$\qed$

\subsection{Embedding of affine $\frak{sl}_2$.}

For $\kappa \in \Cplx$ let  $\cA_\kappa(\frak{sl}_2)$  denote the  
vertex algebroid over $\Cplx$
equal to  $\frak{sl}_2$ as a space,
with bracket $g \zero g' = [g, g']$ and pairing
$g \one g' = \kappa\pair{g | g'}$
where $\pair{ \cdot | \cdot  }$ is the canonically normalized invariant form 
(for $\frak{sl}_2$, it is  $\frac{1}{4} \pair{ \cdot | \cdot  }_{Killing}$ ).

Let 
\begin{equation}
\begin{split}  
\label{sl-2-embedding}
e   & = \de_x  \\
h   & =  - 2 \de_x \opm x  + \la^*    \\
f  &  =  - \de_x \opm x^2  - 2 dx   + x  \la^*   + \half \pair{\la^* | \la^*}  dx
\end{split}
\end{equation}

\begin{lem}
The elements $e$, $f$, $h$ given by the  formulas (\ref{sl-2-embedding})

(1) satisfy the relations of $\widehat{\frak{sl}}_2 (\kappa)$ 
where
 $\kappa = \frac{\pair{\la^* | \la^*}}{2} -2$;

(2) belong to $H^0(\prline,  \cA^{tw}_{\pair{,}} )$
\end{lem}

{\em Proof.}
Restricted to the big cell, 
the statement of Part (1) 
goes back to Wakimoto \cite{W}; see also \cite{F1}.

The rest follows from the following equalities over $U_0 \cap U_1$:
\begin{equation}
\begin{split}  
\de_x   &= - \de_y \opm y^2  -2 dy   + y  \la^*   + \half \pair{\la^* | \la^*}  dy 
\\
 - 2 \de_x \opm x  + \la^* & =   2 \de_y \opm y  - \la^*
 \\
 - \de_x \opm x^2  -2 dx   + x  \la^*   + \half \pair{\la^* | \la^*}  dx  
  & 
  =
   \de_y
\end{split}
\end{equation}
$\qed$

\smallskip

\begin{cor}   
\label{global_sections_of_algebroid-Cor}
The formulas  (\ref{sl-2-embedding}) define an isomorphism  of vertex algebroids over $k$
\begin{equation}
\cA_\kappa (\frak{sl}_2) \iso H^0(\prline, \cA^{tw}_{\pair{, }})
\end{equation}
that extends to the vertex algebra embedding
\begin{equation}
\label{sl-2-vertexalgebra-embedding}
V_\kappa (\frak{sl}_2) \To H^0(\prline, U(\cD^{ch, tw}_{\pair{,}}))
\end{equation}
\end{cor}

{\em Proof.} The map defined by (\ref{sl-2-embedding}) is clearly injective
and the first statement  follows by dimension count. 
The restriction of the second map to the big cell was shown in 
 \cite{F1} to be injective.
 $\qed$

\bigskip


\subsection{The case of a general flag variety}

Recall that we have an identification
\begin{equation}
	\label{Equation isom-hdual-H1}
	\bar{\alpha}:  \fh^* \iso  H^1 (X, \Omega^{1,cl}) \iso {H}^2 (X, \Cplx)
\end{equation}

In other words, the tdo on $G/B$ are classified by $\fh^*$.  
The Lie algebroid $\cT^{tw}_{G/B}$ is an extension
$$
\begin{CD}
0 @>>> \cO_{G/B} \tensor_\Cplx \fh  @>>> \cT^{tw}_{G/B} @>>> \cT_{G/B}  @>>> 0
\end{CD}
$$

A deformation of TCDO is therefore a vertex extension of
$(\cT^{tw}_{G/B}, \pair{,})$ where $\pair{,}$ is a symmetric bilinear pairing
$
\pair{,} : \,  \fh \times \fh  \to \Cplx
$

 We have the following 
\begin{thm}
Let $X = G/B$. 
Then the class of $\cT \CDO^{\pair{}}$
is equal to $0$ if and only if 
$\pair{,}$ is proportional to the restriction of the Killing form on $\fh$.
\end{thm}


{\em Proof.}
First, we find a convenient cocycle representation of the obstruction.

Let $\{ \chi_r  \}$ be the set of fundamental weights,
 $\cL_r$ the corresponding line bundles over $X$,
 $\cD_{\chi_r}$ algebras of tdo acting on $\cL_r$ 
 and
 $T_{\chi_r}$ the  corresponding Lie algebroids. 
 Define the cocycles   
$\mu_r = ( \mu_r^{ij} ) \in \check{Z}^1 (X, \Omega_X^{1, cl})$
 corresponding to $T_{\chi_r}$.
%
  	Then the  map (\ref{Equation isom-hdual-H1}) is the one taking 
	$\chi_r$ to the class of $(\mu_r^{ij})$ in  $ {H}^1 (X, \Omega_X^{1, cl})$.

Take
 $\la^*_r$  to be  the basis of $\fh$ 
dual to the basis $\{ \chi_r  \}$.

  \smallskip
  
  Using Theorem \ref{class-DTCDO-LT-thm}
and the existence of CDO on $X$ (\cite{GMSII}), 
we conclude that  the class of $\cT \CDO^{\pair{}}$ is represented by a cocycle
 $ \lars \, \mu_r^{ij} \wedge \mu^{jk}_s$.
%
%
%
%
%
Its image 
under the natural embedding 
  $H^2(X, \Omega^2 \to \Omega^{3, cl})  \to  H^4(X, \Cplx)$ 
  (cf. \cite{GMSII})
  equals to that of  the element
$$
 S = \lars \, \chi_r \cdot \chi_s   \in  S^2 \fh^*.
$$
which naturally corresponds to the form $\pair{,}: \fh \times \fh \to \Cplx$.

  By \cite{BGG}, $S$ becomes zero in $H^4(X, \Cplx)$  
  if and only if $S$ is $W$-invariant. 
  Therefore, the form $\pair{,}$ has to be a multiple of the Killing form.
  $\qed$

\subsubsection{Embedding of the affine vertex algebra $\cV_{k} (\fg)$}
\
Let $X$ be a $G$-variety.

Let $\ucA_k(\fg)_{X}$ denote the constant sheaf with sections equal to  $\fg$,
equipped with the structure of a $\Cplx_X$-vertex algebroid as follows:
\begin{equation}
\begin{split}
x \zero y   &  = [x, y]  \\
x \one y  &  =  k \pair{x | y} \\
\pi = 0, & \ \ \de=0
\end{split}
\end{equation}
Let $\cA$ be a   (locally trivial) $\pair{,}$-deformation of TCDO.

\smallskip

Let us assume that there is a Lie algebra morphism
\begin{equation}
\alpha: \fg \to \cT^{tw}_X
\end{equation}
lifting the morphism $\bar{\alpha} :\fg \to \cT_{X}$ induced by the action of $G$.
(This is the case for $X = G/B$).

\smallskip

Consider the sheaf of homomorphisms  of vertex algebroids
	\begin{equation}
	 \label{sheaf of embeddings}
	  	\calH om_{\alpha}  ( \ucA_k(\fg)_{X},  \cA)
	\end{equation}
	that lift the morphism $\alpha$.

We are mainly interested in the global sections of this sheaf, as they correspond to embeddings of the vertex algebra $\cV_{k} (\fg)$ 
into the envelope of $\cA$.

\begin{prop}
Suppose the image of $\fg$ in $\cT_X$ generates $\cT_X$ as an $\cO_X$-module.

Then
the sheaf (\ref{sheaf of embeddings}), if locally nonempty,
is an $\Omega^{2,cl}$-torsor.
\end{prop}
	
{\em Proof.}
Let us work locally on a subset $U\subset X$ small enough to admit an identification
$\cA |_U \iso \cT_X^{tw}|_U  \oplus \Omega^1_U$.

Let $w, w' \in \calH om_{\alpha}  ( \ucA_k(\fg)_{X},  \cA) (U)$.

Then $w'(g) = w(g) + \omega(g)$
for some $\omega : \fg \to \Omega^1$, since the $\cT^{tw}$-component is fixed.

Analysis similar to that in \cite{MSV, GMS} shows that 
$\omega$ must be given by 
$$
\omega (g) = \iota_{\alpha(g)} \beta
$$
where $\beta \in \Omega^{2, cl}_X$.

{\sloppypar
Conversely, adding $ \iota_{\alpha( - )} \beta$ to any 
$w \in \calH om_{\alpha}  ( \ucA_k(\fg)_{X},  \cA) (U)$
gives an element of 
$\calH om_{\alpha}  ( \ucA_k(\fg)_{X},  \cA) (U)$.
The statement follows.
$\qed$
}

\smallskip

\begin{rem}
When $\dim X = 1$ the torsor (\ref{sheaf of embeddings}) is trivial, therefore the existence of local embeddings implies 
the existence of a global one. 
For a general flag variety we do not know whether the torsor (\ref{sheaf of embeddings}) 
is trivial, but we believe it is.
\end{rem}



\begin{thebibliography}{99}

	\bibitem[AChM]{AChM} T.~Arakawa, D.~Chebotarov, F.~Malikov,
	Algebras of twisted chiral differential operators and affine localization of $\fg$-modules,
	arXiv: 0810.4964, 
	to appear in Selecta Math.

	
	\bibitem[BB1]{BB1} A.~Beilinson, J.~Bernstein,
	Localisation de $\frak{g}$-modules. (French)  C. R. Acad. Sci. Paris S\'er. I Math.
	292  (1981),  no. 1, 15--18.

	
	\bibitem[BB2]{BB2} A.~Beilinson, J.~Bernstein,
	 A proof of Jantzen conjectures. I. M. Gelfand Seminar, 1--50,
	 Adv. Soviet Math., 16, Part 1, Amer. Math. Soc., Providence, RI, 1993.


	\bibitem[BD1]{BD1} A.~Beilinson, V.~Drinfeld,  Chiral algebras. American
	Mathematical Society Colloquium Publications, 51. American
	Mathematical Society, Providence, RI, 2004. vi+375 pp. ISBN:
	0-8218-3528-9


	\bibitem[BGG]{BGG}
	J.~Berstein, I.M.~Gelfand, S.I.~Gelfand.    
	Shubert cells and the cohomology of the flag manifold (Russian). 
	Func.Analysis and  Applications, vol.7,  no. 1, (1973), 64-65.



	\bibitem[Bre]{Bre} P.~Bressler, The first Pontryagin class. Compositio
	Mathemarica,  {\textbf{143}} (2007), 1127-1163
	
	\bibitem[D]{D} 
	 P.~Deligne, La formule de dualit\'{e} globale, 1973. SGA 4 III, Expos\'{e} XVIII.


	\bibitem[F1]{F1} E.~Frenkel, Wakimoto modules, opers and the center at the
	critical level. Adv. Math. 195 (2005), no. 2, 297--404.

	\bibitem[F2]{F3} E.~Frenkel, Langlands correspondence for loop groups,
	Cambridge University Press, 2007

	\bibitem[FF1]{FF1} B.~Feigin, E.~Frenkel, Representations of affine Kac-Moody
	algebras and bosonization, in: V.Knizhnik Memorial Volume,  L.Brink,
	D.Friedan, A.M.Polyakov (Eds.), 271-316, World Scientific,
	Singapore, 1990

	\bibitem[FF2]{FF2} B.~Feigin, E.~Frenkel, Affine Kac-Moody algebras at the
	critical level and Gelfand-Dikii algebras, in: Infinite Analysis,
	eds. A.Tsuchiya, T.Eguchi, M.Jimbo, {\it Adv. Series in Math. Phys.}
	{\bf 16} 197-215, Singapore, World Scientific, 1992





	\bibitem [GMS1]{GMS}  V.~Gorbounov, F.~Malikov, V.~Schechtman, Gerbes of chiral
	differential operators. II. Vertex algebroids, Invent. Math. 155
	(2004), no. 3, 605-680.

	\bibitem[GMS2]{GMSII} V.~Gorbounov, F.~Malikov, V.~Schechtman,
	 On chiral differential operators over homogeneous spaces.  Int. J. Math. Math. Sci.  26  (2001),  no.2, 83--106.
	

           \bibitem[LWX]{LWX}    Z.-J.~Liu, A.~Weinstein, and P.~Xu,   
           Manin triples for Lie bialgebroids.    
           J. Differential Geom. Volume 45, Number 3 (1997), 547-574. 

	\bibitem[MSV]{MSV} F.~Malikov, V.~Schechtman, A.~Vaintrob,  Comm. in Math. Phys. {\bf 204} (1999), 439-473

	\bibitem[W]{W} M.~Wakimoto, Fock representations of the affine Lie algebra
	 $A\sp {(1)}\sb 1$. Comm. Math. Phys. {\bf 104} (1986), no. 4, 605--609.

\end{thebibliography}
\end{document}